\documentclass[11pt,reqno]{amsart}
\usepackage{amsmath, latexsym, amsfonts, amssymb, amsthm, amscd}
\advance\hoffset -.75cm 
 
\usepackage{amsfonts} 
\usepackage{latexsym} 
\usepackage{amssymb}

\newcommand{\R}{\mathbb R} 
\newcommand{\Z}{\mathbb Z} 
\newcommand{\N}{\mathbb N} 
\newcommand{\Prob}{\mathbb P} 
\newcommand{\E}{\mathbb E} 
\newcommand{\Rd}{{\mathbb R}^d} 
\newcommand{\Zd}{{{\mathbb Z}^d}} 
\newcommand{\bPr}{\mathbf P}
\newcommand{\bE}{\mathbf E}

\newcommand{\eps}{\epsilon} 
\newcommand{\s}{\varphi} 
\newcommand{\La}{\Lambda}
\newcommand{\sg}{\sigma} 
 
\newcommand{\dd}{{\mathrm d}}

\newcommand{\cF}{\mathcal{F}}

\newcommand{\Pone}{\widetilde{\Prob}_N }
\newcommand{\Ptwo}{\Prob }
\newcommand{\Pthree}{\Prob _N}
\newcommand{\Eone}{\widetilde{\E}_N }

\newcommand{\sbs}{\subset}

\newcommand{\ident}{{\mathchoice {\rm 1\mskip-4mu l} {\rm 1\mskip-4mu l}  
{\rm 1\mskip-4.5mu l} {\rm 1\mskip-5mu l}}}

\newtheorem{teo}{Theorem}[section] 
\newtheorem{lem}[teo]{Lemma} 
\newtheorem{cor}[teo]{Corollary} 
\newtheorem{rem}[teo]{Remark} 
\newtheorem{pro}[teo]{Proposition}

\begin{document} 

\title[An interface above a random hard wall] 
{Enhanced interface   repulsion
\\
from quenched hard--wall randomness}
\author[D. Bertacchi]{Daniela Bertacchi}
\address{D. Bertacchi,  Universit\`a di Milano--Bicocca 
Dipartimento di Matematica e Applicazioni, 
Via Bi\-coc\-ca degli Arcimboldi 8, 20126 Milano, Italy 
}
\email{bertacchi\@@matapp.unimib.it}
\author[G. Giacomin]{Giambattista Giacomin}
\address{G. Giacomin, Universit\'e Paris 7 and Laboratoire de
Probabilit\'es et Mod\`eles Al\'eatoires C.N.R.S. UMR 7599, U.F.R.
Math\'ematiques, Case 7012, 2 Place Jussieu, F-75251 Paris, France
\hfill\break
\phantom{br.}{\it Home page:}
{\tt http://felix.proba.jussieu.fr/pageperso/giacomin/GBpage.html}} 
\email{giacomin\@@math.jussieu.fr}
\date{\today}
\begin{abstract}
We consider the {\sl harmonic crystal}, or {\sl massless free field},
$\s=\{\s_x\}_{x\in \Z^d}$, $d\ge 3$, that is the centered Gaussian
field with covariance given by the Green function of the
simple random walk on $\Z^d$.  
Our main aim is to obtain quantitative
information on the repulsion phenomenon 
that arises when we condition $\s_x$ to be larger than
$\sg_x$, $\sg=\{\sg_x\}_{x\in \Z^d}$ is an IID field (which is also
independent of $\s$), 
for every $x$ in a {\sl large} region  $D_N=ND\cap \Z^d$, with $N$ a positive
integer and $D \subset\R^d$  a rather general bounded
subset of $\R^d$. We are mostly motivated by
results for given typical 
realizations of $\sg$ ({\sl quenched} set--up), since the conditioned 
harmonic crystal may be seen as a model for an equilibrium interface,
living 
in a $(d+1)$--dimensional space, constrained not to go below 
 a inhomogeneous substrate that acts as a hard wall. This substrate 
is {\sl mostly flat},
but presents some {\sl rare anomalous spikes}. 
We consider various types of substrate and we observe that
 the
interface is pushed away from the wall {\sl much}
more than in the case of a flat wall as soon as the 
upward tail of $\sg_0$ is heavier than Gaussian, while essentially
no effect is observed if the tail is sub--Gaussian.
In the critical case, that is the one of approximately Gaussian tail,
the 
interplay of the two sources of randomness,
$\s$ and $\sg$, leads to an
enhanced repulsion effect of {\sl additive} type.
This generalizes  work done in the case of
a flat wall
and also in our case the crucial estimates
are optimal Large Deviation type asymptotics as $N\nearrow \infty$
of the probability that $\s$ lies above $\sg$
in $D_N$. We will consider the annealed case too.
It turns out that {\sl quenched} and 
 {\sl annealed} asymptotics coincide and 
this fact plays a role in the proofs and 
concurs to building an understanding of
the phenomenon.
\end{abstract}

\maketitle

\baselineskip .6 cm 

\section{Introduction and main results} 
\label{sec:intro} 
\setcounter{equation}{0}

\subsection{The harmonic crystal}
An harmonic crystal is (for us) 
the  Gaussian random field 
$\s =\{\s _x \}_{x \in \Zd}\in \R^{\Zd}$, $d \ge 3$, 
such that $\E(\s_x)=0$ and $\E(\s_x\s_y)=G(x,y)$ for every
$x$ and $y$ in  $\Zd$, where  $G:\Zd\times \Zd \rightarrow \R^+$
is the Green function of the simple random walk on  $\Zd$:
\begin{equation}
\label{eq:Gcovariance}
G(x,y)=\sum_{n=0}^\infty p_n(x,y),
\end{equation}
where $p_n(x,y)$ is the probability that
the simple random walk $\{X_k\}_{k=0,1,\ldots}$, with $X_0=x$
and hopping to nearest neighbor sites with
probability $1/2d$, is at site $y$ after $n$ time steps. 
We remark that $G(x,y)=G(x-y,0)=G(y-x,0)$
for every $x$, $y\in \Zd$. In short we will
write $\s \sim {\mathcal N}(0,G(\cdot,\cdot))$: the same notation
will be used for (finite dimensional) Gaussian vectors. 
We observe that
$\s$ is a Gibbsian field (cf. \cite{cf:Georgii}) and can
be characterized by its one point conditional probability:
for every $x$ and every measurable bounded function $h:\R \to \R$
\begin{equation}
\label{eq:Gibbs1}
\E \left[
h(\s_x) \vert \mathcal{F}^\s_{\{ x\}^\complement} 
\right]
(\psi)=
E(h(Z_\psi)), \ \ \ \Prob (\dd \psi ) - \text{a.s.},
\end{equation}
where $Z_\psi\sim {\mathcal N}( (1/2d)\sum_{y:\vert y-x \vert=1} \psi_y, 1)$
and, for $A\subset \Zd$, $ \mathcal{F}^\s_{A}=\sigma (\s_x:x\in A)$.
In equation \eqref{eq:Gibbs1} we also introduced the notation
$E$ ($P$) for expectation (probability) of the random variable involved:
if we need to insist on the measure (say $\mu$) on the probability
space we write $E_\mu$ ($P_\mu$) .
We will reserve the use of $\E$ (and $\Prob$)
for the random field $\s$. 
Notice that from \eqref{eq:Gibbs1}
one easily extracts
the fact that $\s$ is also a Markov field.

In more informal way we may simply say that
$\s$ is a Gibbsian field with respect to the 
formal Hamiltonian
\begin{equation}
\label{eq:Hformal}
H(\s)=\frac 1{8d}\sum_{x,y:\vert x-y\vert=1}
\left(\s_x -\s_y\right)^2.
\end{equation}
This imprecise statement helps getting an intuitive 
grasp on the special features of the harmonic crystal.
We stress in particular two  
facts (see \cite{cf:Georgii},
particularly Ch.~13, for a detailed treatment):

\medskip
\begin{itemize}
\item 
Existence of a Gibbsian field associated to a certain 
$H$ is not guaranteed in $\R^{\Z^d}$, due to the lack
of compactness and  general results to tackle this 
problem
do not apply to the case of \eqref{eq:Hformal}.
As a matter of fact there exists no Gibbs measure
associated to \eqref{eq:Hformal} if $d=1,2$.
Of course the fact that $H$ is a quadratic form
allows for a full solution of the existence problem
and the characterization of the space of all the
Gibbs measures (associated to $H$): in particular
one easily arrives to formula \eqref{eq:Gcovariance} for 
the covariances and one understands the necessity
of being on a lattice in which a simple symmetric 
random walk is transient in order to have existence
of a (infinite volume) Gibbs measure.  
 \item The space of Gibbs measures associated to $H$ is extremely large
(as soon as it is non empty, of course).
One can show that $\mu$ is an extremal element
of such (convex) space of  Gibbs measures if and only if $\mu \sim
\mathcal{N}(u_\cdot, G(\cdot,\cdot) )$, with $(\Delta u)_x=
(1/2d)
\sum_{e\in \Z^d: \vert e \vert =1} (u_{x+e}-u_x)=0$ for every $x$,
that is $u$ is harmonic.
In particular we may choose 
$u_x=a+v\cdot x$ for every choice of $a \in \R$ and
$v\in \R^d$. 
We can interpret $\s_x$ as the height of the interface
above a reference plane: to a certain extent 
the {\sl richness} of the Gibbs space is intimately connected
with 
the interest of the model
as a very simpliflied caricature of a physical
 interface (this issue is developed at length
in \cite{cf:Abraham} and  \cite{cf:FFS}).
\end{itemize}
\medskip

As the reader may have noticed, we have made the choice
not to distinguish between random and numerical variables
when talking about $\s$.

\bigskip

\subsection{A model for entropic repulsion: the case
of a flat wall} 
\label{sec:modflat}
In \cite{cf:BDZ} (see however \cite{cf:Bolthausen} for a review
of the various improvements obtained since then) the authors
considered the problem of identifying the asymptotics
of the probability of the event
\begin{equation}
\Omega_N^+=\left\{\s: \s_x \ge 0 \text{ for every } x\in D_N
\right\}
\end{equation}
where $N$ is a positive integer, $D_N= N D \cap \Z^d$, $D=(-1/2,1/2)^d$,
and the asymptotics is with respect to  $N \nearrow \infty$.
Their main results are essentially two, we restate them here
informally:
\medskip
\begin{itemize}
\item In the sense of exponential asymptotics,
$\Prob(\Omega_N^+)$ behaves like $\exp(-\alpha N^{d-2}\log N)$.
The constant $\alpha$ has been determined:
\begin{equation}
\alpha=2G  \mathrm{Cap}(D),   \ \ \ \ \ \ \ 
G(0,0)=G,
\end{equation}
where Cap$(D)$, $D$ an open  subset of $\R^d$
 is the {\sl Newtonian capacity} of $D$:
\begin{equation}
\label{eq:capacity}
\begin{split}
\mathrm{Cap}(D)&=\inf\left\{\frac 1{2d}
\Vert \partial f \Vert _2^2 \,: \,
f\in C^\infty_0(\R^d; [0,\infty)), \ 
f(r)=1 \text{ for all } r \in D\right\},
\\
&=
\sup_{f\in L^\infty(D)}
\frac{\left(\int_D f (r) \dd r\right)^2}{ \int_D \int_D
f(r) f(r^\prime) R_d \vert r-r^\prime \vert^{2-d} \,
\dd r \dd r^\prime}, 
\end{split}
\end{equation}
in which $\partial$ denotes the gradient, $\Vert \cdot \Vert_2$
is the $L^2$--norm of $\cdot$ and 
\begin{equation}
\label{eq:Rd}
R_d=\lim_{x\to \infty} \vert x\vert^{d-2}G(0,x)\in (0,\infty).
\end{equation}
The equivalence between the two definitions of capacity
in \eqref{eq:capacity} can be found for example in \cite[Lemma~A.8]{cf:BD}
and the existence of the non--degenerate limit
in \eqref{eq:Rd} we refer to \cite{cf:Lawler}.
\item
The trajectories of the field $\s$ that are typical with respect 
to $\Prob (\dd \s \vert \Omega_N^+)$ are {\sl pushed}
to infinity as $N \nearrow \infty$ in the sense
that (\cite[Prop.~1.3 and Lemma~4.7]{cf:BDZ} and \cite[Lemma~3.3]{cf:DG})
\begin{equation}
\label{eq:BDZDG}
\lim_{N\to\infty} \sup_{x\in D_N}
\left \vert
\frac{\E \left(\s_x \big\vert
\Omega_N^+ \right)}{ \sqrt{4G\log N}}
-1 \right\vert=0.
\end{equation}
Essentially what happens is that the field stays flat 
and {\sl does not change its structure}
(see \cite[Th.~3.3]{cf:DG}), but it flees the wall:
and it does this to regain its freedom of 
fluctuating
(this effect is indeed called {\sl entropic repulsion}).
\end{itemize}
\medskip

The two issues are intimately connected. In fact 
having a {\sl good guess} for the behavior of the
trajectories of  
$\Prob (\dd \s \vert \Omega_N^+)$ leads to
a good lower bound on 
the asymptotics of $\Prob (\Omega_N^+)$ (and may suggest
a strategy for the upper bound). On the other side 
the same probability asymptotics
enter in a crucial way in proving that the 
{\sl good guess} on 
$\Prob (\dd \s \vert \Omega_N^+)$
is  really 
{\sl close to} $\Prob (\dd \s \vert \Omega_N^+)$
itself. The asymptotic behavior of $\Prob (\Omega_N^+)$ is 
 however not the only
ingredient and this {\sl two way argument} (from probability
estimates to path properties and viceversa) is by no
means general.
For further discussions on physical aspects 
of entropic repulsion we refer for example
to \cite{cf:BeMF}, \cite{cf:LM} and the several references
therein.

\medskip
To understand  the  results in \cite{cf:BDZ}  one may
start with the most naive guess  for the behavior
of $\s_x$ under $\Prob (\dd \s \vert \Omega_N^+)$ 
in the repulsion region: if in the region $D_N$
the field just translates up of $\min_{x\in D_N}
\s_x$, whose typical behavior under $\Prob (\dd \s)$
is approximately $\sqrt{2dG \log N}$,
then the field would not be bothered by the 
presence of the wall. This (sloppy) argument
is easily translated into a rigorous
lower bound on the probability of $\Omega_N^+$
(see \cite[Lemma~2.3]{cf:BDZ}): one would
then verify that it is not the optimal lower bound.
That translating to height $\sqrt{2dG\log N}$
is not a very good guess is also clear from the
result on the typical height of the field 
that we just mentioned above: the field
moves up to (and not beyond!) $\sqrt{4G\log N}$.
Therefore the guess that the field simply
translates globally isn't really correct
and something more complex is happening.
What happens can be synthetised in the following
way: it does not cost too much (in a Large Deviation
sense) to modify 
extrema of the $\s$ field (they happen only on
one site). But we can go beyond:
preventing the field $\s_x$ from fluctuating freely 
on $o(N^{d-2})$ sites (say: sparsely chosen in $D_N$)
is not a substantial modification. This is a non obvious fact,
and it  is essentially a consequence of the fact that a random walk
that leaves from a site $x\in D_N$ reaches with
probability (almost) one $D_N^\complement$ even if 
it we have put $o(N^{d-2})$ traps in $D_N$. 
It turns out that the typical cardinality of $\{x:
\s_x \le -\sqrt{4G\log N} \}$ is about $N^{d-2}$.
One may therefore {\sl believe} that translating
by slightly more than $\sqrt{4G\log N}$ should suffice: 
on (and around) the
{\sl rebel sites} something different happens, but apart
from these sites, that are {\sl few}, we should
still believe that translating is a good guess.
In the present paper we present a proof of the lower bound
that implements in a direct way this heuristics
and that we believe is more direct than the original
proof (and the modified version proposed in \cite{cf:DG}).

\bigskip

\subsection{The model with a random substrate}
In the physical literature much effort is  devoted
 to investigating a variety of random surface phenomena,
including entropic repulsion effects,
 in the presence of a 
{\sl rough} or {\sl disordered} substrate: 
the analysis covers a variety of interface--substrate models
(do for example  
a general {\sl search}
on the physics archive {\it http://xxx.lanl.gov}
for the key--phrase 'disordered substrate'), most
of which seem at the moment out of the reach of  
mathematical treatment.
Here we look for rigorous results
on purely entropic repulsion effects 
 in the presence of 
a disordered  substrate in the simplified 
framework of the high dimensional harmonic crystal.

This substrate, or wall,
 will be  modelled via a random field $\sg=\{ \sigma \}_{x\in \Zd}
\in \R^{\Zd}$: the law of $\sigma$ 
will be denoted by $\bPr$ ($\bE$). 
The hypotheses on $\sg$  are:
\bigskip

\newcounter{Lcount}
\begin{list}{H.\arabic{Lcount}}
      {\usecounter{Lcount}
      \setlength{\rightmargin}{\leftmargin}}
    \item \label{hp:hyp1} {\it Independence}: $\sg$ is an IID field.
\item \label{hp:hyp2} {\it Almost Gaussian behavior of upward  (or $\sg^+$)
 tails}: 
there exists $Q>0$ such that 
\begin{equation}
\lim_{r\to \infty}
\frac1{r^2}
\log \bPr \left(\sg_0> r\right)=-\frac{1}{2Q}.
\end{equation}
\item \label{hp:hyp3} {\it Weak control on downward  (or $\sg^-$) tails}:
$\bE (\sg_0 ^-)<\infty$. 
  \end{list}

\bigskip
  
Examples of $\sigma$ fields of course include the case of
$\sigma_0 \sim {\mathcal N}(0, Q)$, $Q>0$, or  the absolute value 
 of such a variable. We 
 discuss in Subsection~\ref{sec:heur}
each one of this hypotheses. The model
corresponding to H.\ref{hp:hyp2} 
turns out to be
 the most interesting, but for completness we 
consider also the cases:

\medskip

\newcounter{Lcount4}
\begin{list}{H.2--\arabic{Lcount4}}
      {\usecounter{Lcount4}
      \setlength{\rightmargin}{\leftmargin}}
    \item \label{hp:hyp2-1}
{\it Sub--Gaussian behavior of upward  tails}: 
\begin{equation}
\lim_{r\to \infty}
\frac1{r^{2}}
\log \bPr \left(\sg_0> r\right)=-\infty.
\end{equation}
    \item 
\label{hp:hyp2-2} 
{\it Super--Gaussian behavior of upward tails}: 
there exists $\beta\in (0,1)$  and $Q>0$ such that 
\begin{equation}
\lim_{r\to \infty}
\frac1{r^{2\beta}}
\log \bPr \left(\sg_0> r\right)=-\frac{1}{2Q}.
\end{equation}
\end{list}
\medskip

Also for $\sigma$ our notation 
does not distinguish between
random and numerical variables. 

\medskip

The random fields $\s$ and $\sg$ are assumed to be independent
of each other. We may therefore think the configuration
space to be $\R^\Zd \times \R^\Zd$, endowed with the (local)
product topology and equipped with the Borel $\sigma$--algebra:
on this space the measure is $\bPr \otimes \Prob$.  
Therefore $( \sg, \s)\in \R^\Zd\times \R^\Zd$ 
is a wall--interface configuration. 
We 
introduce an interaction between
$\s$ and $\sg$ by conditioning with respect to a suitable event: 
given $\sigma \in \R^{\Zd}$ and  $A\subset \Zd$, the  
$\sigma$--{\sl entropic repulsion event} on $A$ is defined by
\begin{equation} 
\Omega_{A, \sigma}^+= \left\{\s: \s_x \ge \sg_x \text{ for every }
 x \in A
\right\}.
\end{equation}
We 
mostly impose the repulsion on a rather general domain
$D_N=ND \cap \Zd$, $D$ a bounded connected domain with piecewise
smooth boundary and containing the origin:
we use the shortcut 
notation $\Omega^+_{N, \sg} =
\Omega^+_{D_N, \sg}$. 

\medskip
 We 
talk about {\sl quenched results} 
in the cases in which a $\bPr$--typical
configuration $\sigma$ is chosen and kept 
fixed (while $\s$ is considered random):
in this case we prefer to work 
on the measure space $(\R^{\Zd}, \mathcal{B}(\R^{\Zd}),
\Prob)$ rather than introducing complicated conditioning
notations. Of course it is in this
{\sl quenched set--up} that $\Omega_{A, \sigma}^+$ is an event.

We 
talk instead of {\sl annealed results}
when both $\sg$ and $\s$ are averaged at the same time.
In the annealed set--up, with abuse of notation,
 $\Omega_{A, \sigma}^+$ 
is rather
the event $\{(\sg,\s):  \s_x \ge \sg_x \text{ for every }
 x \in A
\}$.

\bigskip
\subsection{Main results: the case of almost Gaussian $\sg^+$ tails}
One of the main results that we are going to prove is that
the quenched 
probability of the entropic repulsion event
is vanishing exponentially and we 
 identify its 
asymptotic 
behavior. Moreover quenched and annealed asymptotics coincide.

\bigskip

\begin{teo}
\label{th:main1} Under hypotheses H.\ref{hp:hyp1},  H.\ref{hp:hyp2} and  H.\ref{hp:hyp3}
we have that 
\begin{equation}
\begin{split}
\label{eq:quenann}
\lim_{N \to \infty}
\frac{1}{N^{d-2}\log N}
\log \Prob \left( \Omega_{N,\sg}^+\right)&=
\lim_{N \to \infty}
\frac{1}{N^{d-2}\log N}
\log \bPr \otimes \Prob \left( \Omega_{N,\sg}^+\right)
\\ &=-2(G+Q) \mathrm{Cap}(D).
\end{split}
\end{equation}
$\bPr (\dd \sg)$--a.s..
\end{teo}
\bigskip

The proof of Theorem \ref{th:main1} requires the sharpest 
probability estimates
that we obtain in this work  and 
sheds light on the behavior of the conditional
measure $\Prob(\cdot \vert \Omega_{N,\sg}^+)$: this is
the measure that contains the information directly
related to the physical situation we are modelling.
The next result concerns the asymptotics of this measure.
For $\epsilon>0$ and a configuration $\s$ call $N_\epsilon (\s)$
the cardinality  of the set 
\begin{equation}
\left\{ x\in D_N :
\left\vert \frac{\s_x}{\sqrt{4(G+Q)\log N}} -1 \right\vert \ge 
\epsilon \right\}
\end{equation}

\bigskip 

\begin{teo}
\label{th:height} 
Under hypotheses H.\ref{hp:hyp1},  H.\ref{hp:hyp2} and  H.\ref{hp:hyp3},
$\bPr(\dd\sigma)$--a.s. 
for every $\epsilon>0$, $N_\epsilon(\s)/\vert D_N\vert$ tends to $0$ in
probability with respect to $\Prob(\dd \s\vert \Omega_{N,\sg}^+)$.
\end{teo}

\bigskip

We refer to  Section~\ref{sec:path} for further results
on $\Prob(\dd \s\vert \Omega_{N,\sg}^+)$.

\subsection{Super/Sub--Gaussian $\sg ^+$ tails}
Call $N_\epsilon ^\kappa (\s)$, $\kappa=0,1$, the cardinality of the set
\begin{equation}
\left\{ x\in D_N :
\left\vert \frac{\s_x}{(4((1-\kappa)G+ \kappa Q)
\log N)^{(\kappa\beta^{-1} 
+ (1-\kappa))/2}} -1
\right\vert \ge 
\epsilon \right\}.
\end{equation}

\bigskip
\begin{teo}
\label{th:mainprime}
Assume H.\ref{hp:hyp1} and H.\ref{hp:hyp3}.
\begin{enumerate}
\item
Under Hypothesis H.2--\ref{hp:hyp2-1},
$\mathbf{P} (\dd \sg)$--a.s. we have that
\begin{equation}
\begin{split}
\label{eq:quenannprime1}
\lim_{N \to \infty}
\frac{1}{N^{d-2}\log N}
\log \Prob \left( \Omega_{N,\sg}^+\right)&=
\lim_{N \to \infty}
\frac{1}{N^{d-2}\log N}
\log \bPr \otimes \Prob \left( \Omega_{N,\sg}^+\right)
\\ &=-2G \mathrm{Cap}(D),
\end{split}
\end{equation}
and for every $\epsilon>0$, $N_\epsilon^0 (\s)/\vert D_N\vert$ tends to $0$ in
probability with respect to $\Prob(\dd \s\vert \Omega_{N,\sg}^+)$.
\item
Under Hypothesis H.2--\ref{hp:hyp2-2},
$\mathbf{P} (\dd \sg)$--a.s.
we have that 
\begin{equation}
\begin{split}
\label{eq:quenannprime2}
\lim_{N \to \infty}
\frac{1}{N^{d-2}(\log N)^{1/\beta}}
\log \Prob \left( \Omega_{N,\sg}^+\right)&=
\lim_{N \to \infty}
\frac{1}{N^{d-2}{\log N}^{1/\beta}}
\log \bPr \otimes \Prob \left( \Omega_{N,\sg}^+\right)
\\ &=-(4Q)^{1/\beta}\frac{\mathrm{Cap}(D)}{2},
\end{split}
\end{equation}
and for every $\epsilon>0$, $N_\epsilon^1 (\s)/\vert D_N\vert$ tends to $0$ in
probability with respect to $\Prob(\dd \s\vert \Omega_{N,\sg}^+)$.
\end{enumerate}
\end{teo}

\bigskip
We see therefore that in the sub--Gaussian regime
the behavior is not far from the one found in the
case of a flat wall (and in fact, under stronger conditions
on the law of $\sg^-_0$, this part of the result
is a direct consequence of the results in \cite{cf:BDZ}, see 
Section~\ref{sec:mainprime}).
But in the super--Gaussian regime the fluctuations
of the substrate are dominating: we can say that in this 
regime the entropic contribution to the phenomenon
is, to leading order, coming from $\sg$, while of course
the energy contribution is still coming from the $\s$--field
and it appears in the capacity term. 

\subsection{On the results, on
 the strategies of proof and possible generalizations}
\label{sec:heur}
First of all we try to extend the heuristic ideas 
that we sketched at the end of Subsection~\ref{sec:modflat}.
We assume H.\ref{hp:hyp2}.
and we start with an observation that seems to
suggest that the effect of a random quenched hard wall
should be the same of  that of a perfectly flat wall:
it is an immediate consequence of the results in
\cite{cf:BDZ} that if
$\vert \sg_x\vert =o(\sqrt{\log N})$ for every $x\in D_N$,
then one obtains Theorem~\ref{th:main1}  and 
Theorem~\ref{th:height} with $Q=0$: that is the 
phenomenology of the flat wall. 
The argument is not totally convincing, because 
under such conditions on $\sg$ we are on a set which is $\bPr$--negligible
(notice that this is not true under H.2--\ref{hp:hyp2-1}).
However one can  show that a typical  $\sg$ is such that for sufficiently
large $N$ 
the cardinality of 
$\{x\in D_N:\vert \sg_x\vert =o(\sqrt{\log N})\}$ is larger
than $\vert D_N \vert (1-\delta_N)$, for any choice of $\{\delta_N\}_N$,
such that $\delta_N\searrow 0$   and $\delta_N N^\epsilon \to \infty$
for every $\epsilon >0$.
So the game is clearly to understand 
if large excursions of $\sg$, that happen on {\sl thin}
sets, affect the $\s$--field:
 quantitatively we observe that, by Hypotheses H.\ref{hp:hyp1}
and H.\ref{hp:hyp2},
there are about $N^{d-(\alpha/2Q)}$ ($\alpha>0$)
sites $x$ on which
$\sg_x $ is approximately $\sqrt{\alpha \log N}$.
Let us accept  that the $\sg$--levels  with $\alpha>4Q$
do not have any effect (recall the discussion at the end of
Subsection~\ref{sec:modflat}): we remain with 
all the levels with $\alpha\in (0,4Q)$.
Higher levels in principle affect the $\s$ field
more seriously, but they are substantially less than
lower levels: and on the other side one can repeat
a similar discussion for the $\s$ field.
It turns out that the {\sl relevant}
$\alpha$ is $2G/\sqrt{G+Q}$ and these levels
{\sl mostly interact}
 with $\s$--downward spikes of height $\approx 2Q\sqrt{\log N}/\sqrt{G+Q}$, 
and to accomodate both $\sg$ and $\s$, the field  $\s$  translates up to 
$\approx \sqrt{4(G+Q)\log N}$. 
Reasoning this way, the appearance of a final result
that depends only on $(G+Q)$ looks quite miraculous. 

\medskip

The quantity $(G+Q)$ appears naturally if we restrict to the case
$\sigma_0 \sim {\mathcal N}(0,Q)$. Then $\s-\sg \sim {\mathcal N}(0, (G+Q
\mathbf{I}) (\cdot,\cdot))$, with $\mathbf{I} (x,x)=1$  and $\mathbf{I} 
(x,y)=0$ if $x\neq
y$. Observe that the long range part of the covariance
is still given by the Green function (and this
is the part responsible for the appearance  
of the capacity): the large excursions 
essentially depend only on the diagonal 
and this justifies the appearance of $G+Q$.
This is of course not a proof, but
it can be turned into a proof:
note that $\s -\sg$ is an FKG field, see the next subsection,
and apply for example the argument in \cite[\S 4]{cf:DG}
for the lower bound;  a proof of the upper
bound is given in 
Section~\ref{sec:UBrsannealed}. 
 But of course in this case we have solved the 
annealed model and quenched probabilities
may be smaller (Corollary~\ref{th:LBann}).
We have therefore transferred the problem to the
{\sl slippery} issue of {\sl quenched$\,\stackrel{?}{=}\,$annealed}.
We take this occasion to stress that probability estimates
can be really viewed as {\sl free energy} estimates:
one can insert the conditioning with respect
to $\Omega_{N,\sg}^+$ directly in the Hamiltonian,
just by adding the site dependent 1--body
potential $V_x(\s_x)=\infty \ident_{(-\infty, \sg_x)}(\s_x)$.

\medskip

Let us now address the issue of the {\sl
necessity} of the hypotheses on $\sg$ and $\s$:
\medskip
\begin{enumerate}
\item Hypothesis H.\ref{hp:hyp1} can be relaxed and the result
extended to  a large class
of mildly correlated fields. However, even in easy cases (like 
$\sg$ a Markov field with exponentially decaying correlations), the extension 
turns out to be heavy. Moreover one should also
observe that the first part of the proof of the probability
upper bound (Proposition~\ref{th:UBsgm}) fails for strongly correlated
fields. As a matter of fact, it is in the class of strongly
correlated $\sg$--fields that we can exhibit examples in which
{\sl quenched$\neq$annealed} (work in progress).  
\item Hypothesis H.\ref{hp:hyp2}
is not necessary to carry out a full analysis, but as we 
argued, it captures models in which the randomness
of wall and interface {\sl act on the same scale}.
Heavier tails (H.2--\ref{hp:hyp2-2}) lead to a predominance 
of the wall randomness. 
Lighter tails (H.2--\ref{hp:hyp2-1}) 
lead to the phenomenology of the flat wall (to leading order, of course).
\item Hypothesis H.\ref{hp:hyp3}
should not be needed at all. We believe that 
imposing, as an extreme case, the hard wall
condition only with positive probability should not change
the phenomena. However having some a priori lower
bound at every site for  $\s$ under the conditioned
field comes really handy.  
\item We have chosen the most elementary harmonic crystal
to simplify the exposition: essentially nothing changes
if we choose $G(\cdot,\cdot)$ to be the Green function
of a more general symmetric translation invariant irreducible random walk
which performs jumps of finite range $k$. However $\s$ in this
case is $k$--steps Markov and the conditioning arguments
become more cumbersome.  Even cases of infinite
range jumps can in principle be tackled: but then 
one has to resort (as in \cite{cf:BDZ}) to hypercontractive
estimates, while here we simply play on conditioning 
(in a way similar to the case  treated in \cite{cf:BDG}).
\end{enumerate}
\medskip

We conclude this discussion by addressing
the question about the optimality of
Theorem~\ref{th:height}. This theorem should be compared
with the result \eqref{eq:BDZDG} obtained in the case of
a flat wall.
 Since the models coincide
for $Q=0$ (at least if one chooses $\sg_0 \sim {\mathcal N}(0,Q)$)
one naturally suspects that Theorem~\ref{th:height}
could be improved. While in principle 
Theorem~\ref{th:height} should be improvable, for general $G$
and $Q$ one certainly cannot prove a result like 
the one we just mentioned for the flat wall.
Observe in fact that if $dQ>(G+Q)$ extrema of $\sg$ field
pierce the interface and therefore no 
uniformity with respect to $x$ is to be expected,
at least as long as we consider  upper bounds:
{\sl local} (or {\sl almost local})
 deformation of the interface over a sparse
lattice of points, the sites of the large
excursions of $\sg$, are necessarily present.

\bigskip

\subsection{Overview of the sections and some further
notation and preliminaries}
In Section~\ref{sec:LBrs} the main result is the quenched
lower bound on the probability of $\Omega_{N,\sg}^+$:
the annealed bound follows then by a standard argument
(that we detail in Corollary~\ref{th:LBann}).
In Section~\ref{sec:UBrsannealed} we take the opposite
route:
the main result is an annealed upper bound, from which 
the quenched upper bound follows.
Therefore the proof of Theorem~\ref{th:main1} follows from
Proposition~\ref{th:LBsgm}, Corollary~\ref{th:LBann},
Proposition~\ref{th:UBsgm} and Corollary~\ref{th:UBquenched}.

In Section~\ref{sec:path}
we present the proof of Theorem~\ref{th:height}
(which is the combination of
 Proposition~\ref{th:UBlem} and Proposition~\ref{th:ubheight}, along
with some other results, see in particular
Remark~\ref{th:hydrolim}). 

In Section~\ref{sec:mainprime}
we 
sketch the proof of Theorem~\ref{th:mainprime}.

\medskip

An important role is played by the FKG (Fortuin--Kasteleyn--Ginibre)
 inequality (or {\sl positive association} property):
if $E$, $F\subset \R^{\Z^d}$ are two increasing events
($E$ is increasing if $\s\in E$ implies that $\s +\psi\in E$
for every $\psi \in [0,\infty)^{\Z^d}$) then $\Prob (E\cap F) 
\ge \Prob (E) \Prob (F)$. Positively correlated Gaussian 
fields satisfy
the FKG inequality: of this fact 
there exist several proofs (see for example \cite{cf:HP}).

\medskip
We conclude with some notations:
 $\lfloor \cdot \rfloor $ denotes the integer
part of the positive real number $\cdot$. Unless otherwise stated,
$o(1)$ is always considered with respect to $N\to \infty$ (and no uniformity 
should be assumed  with respect to other parameters which may be present).
With standard notation we set $\Phi (z)=\int_{-\infty}^z
(1/\sqrt{2\pi})\exp (-z^2/2 ) \dd z$. We keep $\vert \cdot \vert$
to denote the Euclidean norm, or the absolute 
value in the one--dimensional case: if we write $\Vert r \Vert$ we mean
$\max _{i=1,\ldots, d} \vert r_i \vert$, $r \in \R^d$.  
 For $A \subset \Zd$ we 
denote by $\cF^{\sg, \s}_A$ the $\sigma$--algebra
generated by $\sg _x$ and $\s _x$, $x \in A$, 
and  $\cF^{ \sg}_A$ is
the $\sg$--algebra
generated by the $\sg$--variables indexed by $A$, and 
analogous for $\s$. Moreover if $A$ is missing
in this notation, it means $A=\Z^d$.

\section{Probability lower bounds: quenched (and annealed) estimates} 
\label{sec:LBrs} 
\setcounter{equation}{0} 

In this section we work under the hypotheses H.\ref{hp:hyp1} and H.\ref{hp:hyp2}. 
The main result that we 
prove is the following:

\bigskip
\begin{pro}
\label{th:LBsgm}
$\bPr (\dd \sg )$--a.s.
\begin{equation}
\label{eq:LBsgm}
\liminf_{N\rightarrow \infty}
\frac1{N^{d-2}\log N}\log \Prob
\left( \Omega_{N,\sigma}^+ \right)\ge -2(G+Q) \mathrm{Cap}(D). 
\end{equation}
\end{pro}

\bigskip
\noindent
{\it Proof.} For 
every choice of a  large integer parameter $\overline k$ we define the auxiliary 
field $\widetilde \sg$ by setting $\theta=\sqrt{4Q}(1+(1/2\overline k))/\overline k$,
$\widetilde k = \lfloor (\sqrt{2(d+2)Q})/\theta \rfloor $  and
\begin{equation}
\label{eq:aux}
\widetilde \sg _x=
\begin{cases}
\theta \sqrt{\log N} &\text{if } \sigma_x\le \theta \sqrt{\log N}, \\
k\theta \sqrt{\log N}, 
&\text{if } \sigma_x \in ((k-1)\theta \sqrt{\log N},
k\theta \sqrt{\log N}] \text{ for } k =2,3,\ldots, \overline k,\\
 \widetilde k\theta \sqrt{\log N},
&\text{if } \sigma_x \in (\overline k \theta \sqrt{\log N},
\widetilde k\theta \sqrt{\log N}]   \\
\infty &\text{if } \sigma_x > \widetilde k \theta\sqrt{\log N},
\end{cases}
\end{equation} 
and set $L_N(k)=\{x\in D_N:\widetilde \sg _x=k\theta \sqrt{\log N}\}$
for $k\in  \{1,2,\ldots, \overline k , \widetilde k, \infty\}$.

\medskip
Let us now select a good $\sigma$--set.
 Call $N_k $ the 
cardinality of the random set $L_N(k)$.
First we define ${\tt G}_N \in \sigma (\sigma_x: x \in D_N)$ 
as the event specified by 
\begin{equation}
\label{eq:G1}
\left\vert N_k - \bE [N_k] \right\vert \le 
\frac{\bE [N_k]}2, \ \ \  \text{ for } k=2,3, 
\ldots, \overline k, \widetilde k 
\end{equation}
and by
\begin{equation}
\label{eq:G2}
N_\infty=0.
\end{equation}
The good $\sg$--set is 
$({\tt G}_N, \text{ev})=
\bigcup_N\bigcap_{k\ge N}{\tt G}_k$.

\bigskip
\begin{lem}
\label{th:goodset}
$\bPr \big( ({\tt G}_N, \text{ev})\big)=1$.
\end{lem}
\bigskip

\noindent{\it Proof of Lemma \ref{th:goodset}}. 
Set $p_N^k= \bPr (\sigma_0 \in ((k-1)\theta \sqrt{\log N},
 k\theta \sqrt{\log N}])$ and $f_N(k)= N^{d- ( (k-1)^2\theta^2/(2Q))}$
for $k=2,3,\ldots , \overline k$;
$p_N^{\widetilde k}=\bPr (\sigma_0 \in (\overline k\theta \sqrt{\log N},
 \widetilde k\theta \sqrt{\log N}])$
and $f_N(\widetilde k)=
N^{d- ( \overline k^2\theta^2/(2Q))}$.
Then
$\bE [N_k]=\vert D_N \vert p_N^k$ and $\mathrm{var}_\bPr (N_k)=
\vert D_N \vert p_N^k (1-p_N^k)$, therefore
by assumption~H.\ref{hp:hyp2}
for every $\epsilon>0$ 
we have that
\begin{equation}
\label{eq:Nkest}
\lim_{N \to \infty}
N^{-\epsilon}\left(\frac{\bE [N_k]}{ f_N(k)}\vee\frac{ f_N(k)}{\bE [N_k]}\right)=0
\ \ \ \text{ and } \ \ \ 
 \lim_{N\rightarrow \infty}
\frac{\mathrm{var}_\bPr (N_k)}{\bE [N_k]}=1,
\end{equation}
 for $k=2,3,\ldots, \overline k, \widetilde k$.

We use the following inequality due to G.~Bennett (cf. \cite{cf:Bennett})
that says that if $\{ X_j\}_{j=1,2,\ldots}$ 
is a collection of centered IID variables such that  
$\vert X_1 \vert \le 1$, then  for every $t\ge 0$
\begin{equation}
\label{eq:Bennett}
P\left(
\left\vert \sum_{j=1}^n X_j \right\vert >t \right)
\le 2 \exp
\left(
-\frac {t^2} {2n \mathrm{var}(X_1)+2t/3 }
\right).
\end{equation}
Therefore for $k=2,3,\ldots, \overline k, \widetilde k$
\begin{equation}
\bPr\left(
\left \vert N_k - \bE [N_k] \right\vert 
> \frac{\bE[N_k]}2
\right)\le
2\exp \left(
-\frac {(\bE[N_k])^2}
{
8 \mathrm{var}_{\bPr} (N_k)+{4\bE[N_k]}/3
}
\right),
\end{equation}
and applying \eqref{eq:Nkest}
we obtain that for every sufficiently large $\overline k$
there exists $c>0$ such that for $k=2,3,\ldots, \overline k, \widetilde k$
\begin{equation}
\bPr\left(
\left \vert N_k - \bE [N_k] \right\vert > \frac{\bE[N_k]}2
\right)\le
2 \exp(-c \bE [N_k])\le 2 \exp\left( 
-c N^{(d-2)/2}
\right).
\end{equation}
Moreover by direct computation 
$\bPr (N_\infty>0)\le N^{-3/2}$ for sufficiently large $N$.
The first Lemma of Borel--Cantelli completes the proof.
\qed (Lemma \ref{th:goodset})

\bigskip

From now on we simply assume
that $\sigma \in ({\tt G}_N, ev)$.
So, in particular, $L_N (\infty)=\emptyset$ and 
for every $\epsilon>0$ there exists $\overline{N}$ such
that for $N\ge \overline{N}$ 
\begin{equation}
\begin{split}
N_k & =\vert L_N (k)\vert \le  N^{d- ( (k-1)^2\theta^2/(2Q))+\epsilon},\qquad
k=1,2,\ldots, \overline k,\\
N_{\widetilde k} & =\vert L_N (\widetilde k)\vert \le  
N^{d- (\overline k^2\theta^2/(2Q))+\epsilon}, \\
\end{split}
\end{equation}
(notice that the result is trivial for $k=1$).

Let us go back to the analysis of the $\s$ field:
we have that 
\begin{equation}
\label{eq:sepmaj}
\begin{split}
\Prob \left( \s_x \ge \sg _x, \
x \in D_N \right)
& \ge 
\Prob \left(
\s_x \ge \widetilde \sigma _x, \  x \in D_N
\right)
\\
&\ge 
\Prob \left(   \s_x \ge \widetilde \sg _x ,
\ x \in D_N^- \right)
\cdot\Prob \left(\widetilde \Omega_N^+ (\widetilde k)  \right),
\end{split}
\end{equation}
in which the first step is immediate consequence of $\widetilde{\sg}\ge \sg$
and in the second one 
we used the FKG inequality  with the notations 
$\widetilde \Omega_N^+(\widetilde k) :=\{\s :\s_x\ge \widetilde{k}
\theta\sqrt{\log N},  x\in L_N(\widetilde k)\}$ and  $D_N^-=D_N\setminus 
L_N(\widetilde k)$. 
Therefore 
\begin{equation}
\label{eq:LBsgm1st}
\begin{split}
\liminf_{N\rightarrow \infty}
\frac1{N^{d-2}\log N}  \log \Prob
\left( \Omega_{N,\sg}^+ \right)
\ge &
 \liminf_{N\rightarrow \infty}
\frac1{N^{d-2}\log N} 
\log\Prob \left(\widetilde \Omega_N^+ (\widetilde k) \right)
\\
&+ \liminf_{N\rightarrow \infty}\frac1{N^{d-2}\log N}
\log\Prob \left( \s_x \ge \widetilde \sg _x ,
\ x \in D_N^-
  \right)
.
\end{split}
\end{equation}

The following straightforward entropy estimate
deals with the first term in the right--hand side
of the above expression:
let us introduce 
the map $T_\psi$: $(T_\psi \s )_x= \s_x +\psi_x$, for $\psi \in \Omega_N$ and
$x \in \Z^d$. 
If $\mu$ and $\nu$ are two probability measures defined
on the same measurable  space and if $\mu$ is absolutely continuous
with respect to $\nu$ we 
 denote by $H(\mu \vert \nu)$ the relative entropy 
$E_{\mu} [\log (\dd \mu /\dd \nu)]$.
If we choose $\psi_x=  \sqrt{2G(d+2) \log N} +
\widetilde k \theta \sqrt{ \log N}$ for
$x\in L_N(\widetilde k)$ and $\psi_x=0$ otherwise, 
by direct computation
\begin{equation}
H\left( \Prob T_\psi^{-1} \vert \Prob 
\right)= \frac 1{4d} \sum_{x,y: \vert x-y \vert=1}
\left( \psi_x -\psi_y \right)^2\le \left(  \sqrt{2G(d+2) \log N} +
\widetilde k \theta \sqrt{ \log N}    \right)^2
 \vert L_N (\widetilde k) \vert.
\end{equation}
By \eqref{eq:G1} and \eqref{eq:Nkest} we therefore have that
for $N$  and $\overline k$ sufficiently large
\begin{equation}
H\left( \Prob T_\psi^{-1} \vert \Prob 
\right)\le N^{d-2(1+ (1/3\overline k))^2}.
\end{equation}
Moreover by using the FKG inequality we obtain that
for sufficiently large $N$ 
\begin{equation}
\begin{split}
\Prob T_\psi^{-1}
\left(\widetilde \Omega_N^+ (\widetilde k)  \right)&=
\Prob\left( \s_x \ge -  \sqrt{2G(d+2) \log N}
\text{ for every } x \in L_N (\widetilde k)\right)
\\
&\ge \prod_{ x\in L_N (\widetilde k)}
\Prob\left( 
\s_x \ge - \sqrt{2G(d+2) \log N}
\right)\ge ( 1-(1/N^{d+1}))^{N^d}\ge 1/2,
\end{split}
\end{equation} 
and therefore by applying the standard entropy inequality
\begin{equation}
\label{eq:entineq}
\log \left(
\frac{P_\nu (E)}{P_\mu (E)}
\right)\ge-\frac{ 1} { P_\mu (E)}
\left( H(\mu\vert \nu) +e^{-1}\right),
\end{equation}
in which $\mu$ and $\nu$ are two probabilities and 
$E$ is an event of positive $\mu$ measure,
we obtain that
\begin{equation}
\Prob \left(\widetilde \Omega_N^+ (\widetilde k)  \right)\ge
\exp \left(-N^{d-2(1+ (1/4\overline k))^2}\right),
\end{equation}
for sufficiently large $N$, which 
shows that the first term in the right--hand side of
\eqref{eq:LBsgm1st} vanishes. 
\medskip

Let us therefore concentrate on the second term
and on the event
$\widetilde{\Omega_N^+}=\left\{ \s_x \ge \widetilde \sg _x ,
\ x \in D_N^-
  \right\}$:
the proof of Proposition \ref{th:LBsgm}
is  complete once we have shown that
\begin{equation}
\label{eq:splitlevel}
\liminf_{N \to \infty}
\frac 1{N^{d-2}\log N}
\log \Prob\left(\widetilde{\Omega_N^+}\right)
\ge -2 (G+Q)\mathrm{Cap}(D).
\end{equation}

\medskip
In order to prove \eqref{eq:splitlevel}  we set 
 $\alpha_N= \alpha \sqrt{ \log N}$, $\alpha >0$. 
For $\psi_N \in \Omega$ such that  $({\psi_N})_x= \alpha_N f(x/N)$,
$f\in C^\infty_0 (\R^d; [0,\infty))$ and $f(r)=1$
if $r \in D$,
we set $\Pthree=\Prob T_{\psi_N} ^{-1}$ and $\Pone(\cdot)=
\Prob_N ( \cdot \vert \widetilde{\Omega_N^+})$.
Therefore 
${\dd\Pone}/{\dd\Ptwo}=(
{\dd \Pone}/{\dd\Pthree})(
{\dd \Pthree}/{\dd\Ptwo})$
and by the entropy inequality \eqref{eq:entineq}
\begin{equation}
\begin{split}
\log\Prob
\left(
\widetilde{\Omega_N^+}\right) & \ge
-H\left(\Pone \vert \Ptwo \right)-e^{-1}\\
&=-H\left(\Pone \vert \Pthree \right)
-\Eone \left( \log \left(\frac{\dd \Pthree}{\dd\Ptwo} \right) \right)-e^{-1}
\equiv -H_1-H_2-e^{-1}.
\end{split}
\end{equation}

First of all by direct evaluation and FKG we have
\begin{equation}
\begin{split}
\label{eq:H1}
H_1=-\log \Prob T_{\psi_N} ^{-1}\left(\widetilde{\Omega ^+_{N}} \right)
&
\le -\sum_{x \in D_N ^-}
\log \Prob \left( \s_x \ge 
\widetilde\sg_x -\alpha_N \right)
\\ 
&= -\sum_{k=1}^{\overline k}\vert L_N(k) \vert \log 
\left( 1-
 \Prob \left( \s_0 < k\theta\sqrt{\log N}  -\alpha_N\right)\right).
\end{split}
\end{equation}
One checks directly  that if
\begin{equation}
\label{eq:alphacond}
\alpha>\overline k\theta ,
\end{equation}
and
\begin{equation}
\label{eq:combcond}
\frac{(k-1)^2 \theta^2}{2Q}
+
\frac{(k\theta-\alpha)^2}{2G}
>2,
\end{equation}
for all $k\le\overline k$, then each of the $\overline k$
summands in \eqref{eq:H1} is $o(N^{d-2})$, and therefore negligible:
\begin{equation}
\label{eq:H1res}
\lim_{N \rightarrow \infty}
\frac1{N^{d-2}}H_1=0.
\end{equation}
Observe that, by \eqref{eq:alphacond} and \eqref{eq:combcond},
a more explicit assumption that implies 
\eqref{eq:H1res} is 
\begin{equation}
\label{eq:condalpha}
\alpha > k\theta  +\sqrt{4G- (k-1)^2\theta^2 \left(
\frac{G}{Q}\right)}, \ \text{ for every } k\le\overline k.
\end{equation}
Note that $4G- (k-1)^2\theta^2 (G/Q)\ge0$ 
holds for every $k\le\overline k$ using the definition of
$\theta$.
If we observe that
\begin{equation}
 \max_{1\le k\le\overline k}
k\theta +\sqrt{4G- (k-1)^2\theta^2 \left(
\frac{G}{Q}\right)}\le \theta + \max_{x \in [0,2\sqrt{Q}]}
\left\{
x+ \sqrt{4G-x^2 \left(\frac{G}{Q}\right)}
\right\}\le \theta +2 \sqrt{G+Q},
\end{equation}
we conclude that 
\eqref{eq:condalpha} is satisfied if
\begin{equation}
\label{eq:condalpha2}
\alpha>2\sqrt{G+Q}+\theta.
\end{equation}
Therefore under this hypothesis on $\alpha$
the estimate \eqref{eq:H1res} holds.

\medskip

Let us consider $H_2$: observe that
\begin{multline}
\log \left(\frac{\dd \Pthree}{\dd\Ptwo} (T_{\psi_N}\s) \right)
=
\\
\frac 1{4d}
\sum_{x,y: \vert x-y \vert=1}
\left(  2\alpha _N (f(y/N)-f(x/N))(\s_x-\s_y)+
\left(\alpha _N (f(x/N)-f(y/N)) \right)^2
\right),
\end{multline}
and therefore
\begin{equation}
\begin{split}
\frac{1}{N^{d-2}\log N}\,
H_2=&
\frac{\alpha_N^2}{4dN^{d}\log N}
\sum_{x,y:\vert x-y \vert=1}
\left(N(f(x/N)-f(y/N)) \right)^2
\\
&+\frac{\alpha _N}{ 2dN^{d-2}\log N}\,
\E\left[
\sum_{x,y:\vert x-y  \vert=1}
(f(y/N)-f(x/N))(\s_x-\s_y)
\bigg\vert 
T_{\psi_N}^{-1}\widetilde{\Omega_N^+}
 \right]\\
\equiv & \, C_N+R_N.
\end{split}
\end{equation}
It is easy to see that $C_N$ converges
for $N \rightarrow \infty$ to 
$\alpha^2 \Vert \partial f \Vert_2^2/(4d)$. 
We show now  that
$\lim_{N \rightarrow \infty} R_N=0$ if 
\eqref{eq:condalpha2} holds.
Observe in fact that
\begin{equation}
\frac 1{2d}\sum_{x,y:\vert x-y \vert =1}
(f(x/N)-f(y/N))(\s_x-\s_y)
=
-\sum_x (\Delta f(\cdot/N))(x) \s _x
\sim {\mathcal N}\left( 0, \sigma_N^2\right),
\end{equation}
where $\sigma_N^2= (\Vert \partial f\Vert^2_2/2d) N^{d-2}
(1+o(1))$.
We use now the following consequence of 
Jensen inequality ($Y $ a random variable,
$E$ a positive probability event, $t>0$)
\begin{equation}
\label{eq:entineq2}
\E\left[
Y \vert E\right]
\le \frac 1t\log \E \left[\exp(tY)\right]
- \frac 1t \log \Prob \left( E \right),
\end{equation}
to obtain with $t=N^{d-2}$ (recall \eqref{eq:H1}) that 
\begin{equation}
\begin{split}
\vert R_N \vert &\le  \frac 1{t}\log \E \left[ 
\exp \left( \frac{t\alpha _N}{ 2dN^{d-2}\log N}
\sum_{x,y}
(f(y/N)-f(x/N))(\s_x-\s_y)
\right)
\right]
- \frac 1t\log \Prob \left( T_{\psi_N}^{-1} \Omega_N^+ \right)
\\
& \le   \frac{ \alpha^2 \sigma_N^2 }{2(2d)^2N^{d-2}\log N}+ 
\frac{H_1}{N^{d-2}}=o(1).
\end{split}
\end{equation}
 This shows that under the hypothesis \eqref{eq:condalpha2} on $\alpha$
 \begin{equation}
 \lim_{N\rightarrow \infty}\frac 1{N^{d-2}\log N}
H_2
=\frac{\alpha^2}2 \frac{\Vert \partial f \Vert_2^2}{2d}.
\end{equation}
Since $\Pone (\Omega_N^+)=1$ 
we may apply the entropy inequality 
\eqref{eq:entineq}, and \eqref{eq:splitlevel}
is obtained by optimising the choices
of $f$ and $\alpha$, by the definition 
of the capacity (cf. \eqref{eq:capacity}, first line) and using the fact that $\theta$ can be taken
arbitrarily small (that is, $\overline k$ arbitrarily large).
\qed

\bigskip

\begin{rem}
One may wonder if a more general estimate
like the one proven in \cite[\S 2]{cf:BDZ} holds in this case too.
The answer is positive: 
Proposition~\ref{th:LBsgm}
can be extended in the sense that
if $\{ b_N\} _{N \in \N}$ is a sequence of real numbers such that 
$\lim_{N\rightarrow \infty} b_N /\sqrt{N}=b\ge -2\sqrt{G+Q}$,
then $\bPr (\dd \sg )$--a.s.
\begin{equation}
\label{eq:LBsgmgen}
\liminf_{N\rightarrow \infty}
\frac1{N^{d-2}\log N}\log \Prob
\left( \s_x \ge \sg _x+ b_N \text{ for every }
x \in D_N \right)\ge -\left(2\sqrt{G+Q} + b\right)^2\mathrm{Cap}(D)/2. 
\end{equation}
We 
do not give a proof of this result, except
for the extreme case $b= -2\sqrt{G+Q}$, which
 follows immediately from 
Lemma~\ref{th:technheight} below.
\end{rem}
\bigskip

In Section~\ref{sec:path}, the proof of an upper bound
on the height of the conditioned field 
requires the following
technical estimate.

\begin{lem}
\label{th:technheight}
For every $\epsilon>0$, if we choose
 $\alpha_N=  (2\sqrt{G+Q}+\epsilon)\sqrt{\log N}$ we have that
\begin{equation}
\label{eq:technheight}
\lim_{N\to \infty}
\frac1{N^{d-2}}
\log
\Prob \left(\Omega_{N,\sigma-\alpha_N}^+\right)=0.
\end{equation}
\end{lem}

\bigskip

\noindent
{\it Proof.}
It is a simplified version of the preceeding proof:
let us keep the same notations.
As for \eqref{eq:LBsgm1st} we have
\begin{equation}
\label{eq:LBsgm1st2}
\begin{split}
\liminf_{N\rightarrow \infty}
\frac1{N^{d-2}}  \log \Prob
\left( \Omega_{N,\sg -\alpha_N}^+ \right)
\ge
&
 \liminf_{N\rightarrow \infty}\frac1{N^{d-2}}
\log\Prob \left( \s_x \ge \widetilde \sg _x -\alpha_N ,
\ x \in D_N^-
  \right)
\\
+ 
 \liminf_{N\rightarrow \infty} &
\frac1{N^{d-2}}  
\log\Prob \left( 
\left\{\s_x\ge \widetilde{k}  \theta\sqrt{\log N}-\alpha_N, 
x\in L_N(\widetilde k)\right\}
 \right)
.
\end{split}
\end{equation}
The second term on the right--hand side of 
\eqref{eq:LBsgm1st2} is not smaller than the first term 
 in  the right--hand side of \eqref{eq:LBsgm1st} and it is therefore
equal to zero. The first term is dealt by applying
the FKG inequality very much in the spirit of
\eqref{eq:H1}:
\begin{multline}
\label{eq:H1P}
\frac1{N^{d-2}}
\log\Prob \left( \s_x \ge \widetilde \sg _x -\alpha_N ,
\ x \in D_N^-
  \right)
\\
\ge \, 
\sum_{k=1}^{\overline k}\vert L_N(k) \vert N^{2-d}
\log 
\left( 1-
 \Prob \left( \s_0 < k\theta\sqrt{\log N}  -\alpha_N\right)\right).
\end{multline}
As for \eqref{eq:H1}, the term on the right--hand side
of  \eqref{eq:H1P} vanishes as $N\to \infty$ if 
\begin{equation}
\frac{(k-1)^2 \theta^2}{2Q}
+
\frac{(k\theta-2\sqrt{G+Q}-\epsilon)^2}{2G}
>2,
\end{equation}
for every $k\le \overline k$.
But this is true as long as $\epsilon>\theta$:
since $\theta$ can be chosen arbitrarily
small, we are done.
\qed

\bigskip

We conclude this section by observing that the quenched
lower bound provides also an annealed lower bound.

\bigskip
\begin{cor}
\label{th:LBann}
\begin{equation}
\label{eq:LBann}
\liminf_{N\rightarrow \infty}
\frac1{N^{d-2}\log N}\log \bPr \otimes \Prob
\left( \Omega_{N,\sigma}^+ \right)\ge -2(G+Q) \mathrm{Cap}(D). 
\end{equation}
\end{cor}

\bigskip
\noindent
{\it Proof.}
Since by Proposition~\ref{th:LBsgm}
for every $\eps>0$ we have that the $\mathbf{P}$ probability that 
$\Prob (\Omega_{N,\sigma}^+)\ge \exp(-2(G+Q+\eps)Cap(D) N^{d-2}\log N)$
is larger that $1/2$ 
for sufficiently large $N$, the result is immediate.
\qed

\bigskip

\section{Probability upper bounds: annealed (and quenched) estimates} 
\label{sec:UBrsannealed} 
\setcounter{equation}{0}

In this section we 
need also some control on the
downward tails of the $\sg$ field. We recall that
in this section we 
commit abuse  of notation for
$\Omega^+_{\Lambda, \sg}$.

\bigskip 
\begin{pro} Under hypotheses H.\ref{hp:hyp1}, H.\ref{hp:hyp2}
and H.\ref{hp:hyp3} we have that
\label{th:UBsgm}
\begin{equation}
\label{eq:UBsgm}
\limsup_{N\rightarrow \infty}
\frac1{N^{d-2}\log N}\log \bPr \otimes\Prob
\left( \Omega^+_{N, \sg} \right)\le  -2\left( G+Q\right)\mathrm{Cap}(D). 
\end{equation}
\end{pro}
\bigskip

\noindent{\it Proof.}
Let us choose an even natural number $L$ and for $y\in 
2L \Zd$ let us set 
\begin{equation}
B(y)=B_{L}(y)=
\left\{
x:\max_{i=1, \ldots , d} \vert x_i - y_i \vert =L/2 
\right\}, 
\end{equation}
and $ \Lambda _c$ is the set of $y\in 
2L \Zd$ such that $B(y) \subset D_N$.
Set also $\Lambda = \bigcup _{y \in \Lambda_c}
B(y)$.
We have that
\begin{equation}
\label{eq:Markov+ind}
\bPr \otimes \Prob 
\left( \Omega^+_{N, \sg} \right)
\le
\bPr \otimes \Prob 
\left( \Omega^+_{\Lambda \cup \Lambda_c, \sg} \right)
=
\bE \otimes \E 
\left[
\prod_{y \in \Lambda _c}
\bPr \otimes \Prob
\left( \s_y  \ge \sg _y \big \vert
\cF ^{\sg, \s}_{B(y)}\right); \Omega^+_{\Lambda , \sg}  
\right],
\end{equation}
in which we used the Markov property of the $\s$--field
and the independence of the $\sigma$--field.
Observe now that, under 
$\bPr \otimes \Prob
(\, \cdot \, \vert   
\cF ^{\sg, \s}_{B(y)})(\psi)$,
$\s_y \sim {\mathcal N}(\sum _{z \in B(y)} q(z) \psi _z , G_L )$
where $q(z)=q_L(z)$ is the
probability that a simple random walk 
leaving at $y$ hits $B(y)$ at $z$
and $G_L$ is a positive number with the property
that $G_L \nearrow G$ as $L\nearrow \infty$. 
We set $M_y^\Box (\psi)= \sum_{z \in B(y)} q(z) \psi _z$.

We now take a positive number $\kappa$  
and we consider the {\sl inner $\kappa$--discretization}  
of $D$: that is
for $r\in \kappa \Zd$, set $A_r=r+[0,\kappa)^d$ and 
define $I=\{r\in \kappa \Zd: 
A_r  \subset D\}$ (assume 
$I \neq \emptyset$).  
We are interested in this decomposition at the lattice 
level or, more precisely, on the $2L$--rarified lattice level (the sublattice 
$\Lambda_c$ of
centers): 
so define $C_r= N A_r \cap \Lambda_c$ and remark that
$\vert C_r \vert=c(N\kappa /2L)^d (1+o(1))$. 
For $\eta\in (0,1/4)$ and  $\alpha\in (0,4(G_L+Q))$ let us now consider the event 
\begin{equation} 
\label{eq:Eea}
E_{\eta, \alpha}= 
\left\{(\sg, \s): \text{there exists }r \in I \text{ such that } 
\vert \{ y \in C_r: M_y^\Box (\s )\le \sqrt{\alpha\log N}\}\vert 
\ge \eta \vert C_r \vert  \right\}. 
\end{equation} 
Of course $E_{\eta, \alpha }$ is $\mathcal{F}^\s$-measurable.
Observe that on $E_{\eta, \alpha }$  
\begin{multline}
\label{eq:twofc} 
\prod_{y \in  \Lambda_c}  
\bPr \otimes \Prob \left(\s_y \ge \sigma_y\big \vert \cF^{\sg, \s } _{B(y)}\right) 
\\
\le \prod_{y \in \Lambda _c} 
\left( 1-
\Prob\left(\s_y\le -\frac{G_L}{G_L+Q}
\sqrt{\alpha \log N} \right)
\bPr \left( \sg_y \ge \frac{Q}{G_L+Q}\sqrt{\alpha \log N}\right)
\right) 
\\
\le
\left( 1-
N^{-\frac{\alpha+\eps}{2(G_L+Q)}
}\right)^{\eta\vert C_{r}\vert},
\end{multline} 
where $r$ is any element of $I$ and $\eps \in (0, 4(G_L+Q)-\alpha)$.
Then for sufficiently large  $N$ and a suitable 
constant $c^\prime$  we have that
\begin{equation} 
\label{eq:neglect1}
\begin{array}{ll}
\bE \otimes \E\left[
\prod_{y \in  \La_c} 
\bPr \otimes \Prob \left(\s_y \ge \sigma_y \vert \cF _{B(y)}\right)
; E_{\eta, \alpha} \right]  &
\le \left( 1-
N^{-\frac{\alpha+\eps}{2(G_L+Q)}}
\right)^{c\eta (N\kappa/L)^{d}}\\
&\le\exp\left(
-{c^\prime} 
N^{d- \frac{\alpha+\eps} {2(G_L+Q)}} \right),\\
\end{array}
\end{equation} 
which is  negligible (recall \eqref{eq:UBsgm} and the choice of $\eps$ and $\alpha$).

Let us consider now the event
\begin{equation}
E_{N}=
\left\{ (\sg, \s):
\left\vert \left\{
y \in \Lambda_c: M_y^\Box (\sg)\le -(\log N)^{1/4}
\right\}\right\vert \ge \delta_N
\vert \Lambda_c \vert \right\},
\end{equation}
with $\delta_N=\max(\sqrt{\bPr(  M_0^\Box (\sg)\le -(\log N)^{1/4})}, N^{-1}
))$.
By H.\ref{hp:hyp3} and the Markov inequality,
$\delta_N$ vanishes as $N\to \infty$.
A direct application of \eqref{eq:Bennett}
leads to the existence of a constant $c>0$ such that
for every $N$
\begin{equation}
\label{eq:neglect2}
\bPr \otimes \Prob
\left( E_{N}\right) \le c\exp(-\delta_N \vert \Lambda_c\vert).
\end{equation}
which again is negligible, in view of the result
we are after (cf. \eqref{eq:UBsgm}).
We observe that on $E_N^\complement$ we may select a set 
(depending on $\sigma$ and $N$)
$\Lambda_c^G \subset \Lambda_c$,
with the property that $\vert \Lambda^G_c \vert / \vert \Lambda_c \vert \ge (1-\delta_N)$,
on which  $M_y^\Box (\sg)> -(\log N)^{1/4}$: this implies that if we define
$C_r^G=C_r \cap \Lambda_c^G$ we can find a positive constant $c=c(\kappa, D)$
such that $\vert C_r^G\vert / \vert C_r\vert \ge (1-c\delta_N)\vee 0$ for every $r$.
We choose $N$ such that $c\delta_N<\eta$.

\medskip

By the last two observations (\eqref{eq:neglect1} and \eqref{eq:neglect2}), 
we are allowed to  replace the event $\Omega^+_{\Lambda , \sg}$
with $\Omega^+_{\Lambda , \sg}\cap E_{\eta, \alpha}^\complement
\cap E_{N}^\complement$
in the rightmost expression in 
\eqref{eq:Markov+ind}.
If $(\sg,\s) \in 
\Omega^+_{\Lambda , \sg}\cap E_{\eta, \alpha}^\complement
\cap E_{N}^\complement$
then
 for every $r \in I$ there are at least 
$(1-2\eta) \vert C_r^G\vert$ sites $y \in C_r^G$ such that 
$M_y^\Box (\s)> \sqrt{\alpha\log N}$ 
and in the remaining  (at most  $2\eta \vert C_r^G\vert$) sites
$M_y^\Box (\s -\sg)\ge 0$ and $M_y^\Box (\sg)> -(\log N)^{1/4}$.
Therefore 
for every choice of $f_r\ge 0$, $r \in I$, 
\begin{equation} 
\label{eq:mc} 
\sum_{r \in I} 
f_r \frac{1}{\vert C^G_r \vert} \sum_{y \in C^G_r} 
M_y^\Box (\s) \ge (1-3\eta) \sqrt{\alpha \log N} 
\sum_{r \in I} 
f_r.
\end{equation} 
Therefore if we call $F_N$ the event
specified by $(\sg, \s)$ such that \eqref{eq:mc} holds, we have 
shown that 
\begin{equation}
\label{eq:UBviaFandE}
\limsup_{N\to \infty}\frac1{N^{d-2}\log N}
\bPr \otimes \Prob \left(\Omega_{N,\sg}^+\right) \le
\limsup_{N\to \infty}\frac1{N^{d-2}\log N}  \log \bPr \otimes \Prob \left(
F_N \cap E_N^\complement \right).
\end{equation}
In order to deal with Gaussian computations we condition
with respect to $\sg$: notice that $E_N$ is measurable with
respect to $\mathcal{F}^\sg$.
We have that on $E_N^\complement$
\begin{equation}
\label{eq:boundF}
\bPr \otimes 
\Prob \left( F_N\big \vert \mathcal{F}^\sg\right)  \le 
\exp \left(
-\frac{(1-3\eta)^2 \alpha\log N \left(
\sum_{r \in I} f_r
\right)^2}{2\mathrm{var}\left(\sum_{r \in I} 
f_r \frac{1}{\vert C^G_r \vert} \sum_{y \in C^G_r} 
M_y^\Box (\s)\right)}
\right).
\end{equation} 
We estimate the variance with respect
to the Gaussian measure $\bPr \otimes 
\Prob \left( \dd \s \vert \mathcal{F}^\sg\right)$: by Jensen's inequality
\begin{equation}
\mathrm{var} \left(\sum_{r \in I} 
f_r \frac{1}{\vert C_r^G \vert} \sum_{y \in C_r^G} 
M_y^\Box (\s)\right) \le
\mathrm{var}_{\Prob}\left(\sum_{r \in I} 
f_r \frac{1}{\vert C_r^G \vert} \sum_{y \in C_r^G} 
\s_y \right),
\end{equation} 
 and observe that, 
 if we define the function $f_\kappa: \Rd \rightarrow \R$
 as $f_\kappa (x)=
 \sum_{r \in I} f_{r}
 \ident_{A_{r}}(x)$,
we can write 
\begin{equation} 
\sum_{r \in I} 
f_r \frac{1}{\vert C_r^G \vert} \sum_{y \in C_r^G} 
\s_y  =
\frac1{\vert C_r \vert}
\sum_{y\in \La^G_c} \gamma_N (y) f_\kappa (y/N) \s _y ,
\end{equation}
where $\gamma_N (y)= \vert C_r \vert/ \vert C^G_{r(y)} \vert$,
with $r(y)$ the index of the $C$--box that contains $y$.
Note that $\gamma_N (y)$ is a $\sg$--dependent function,
but, uniformly in $\sg$, $1\le \gamma_N (y) \le 1/(1-c\delta_N)$
for every $r$ and every sufficiently large $N$ (chosen before).
Observe also that
\begin{equation}
\sum_{r \in I} 
f_r= \frac1{\vert C_r \vert} \sum_{y\in \La_c }  f_\kappa (y/N).
\end{equation}
By direct computation we obtain that
\begin{equation}
\label{eq:dil_cap}
\lim_{N\rightarrow \infty}
\frac1{N^{d-2}}
\frac
{\left(\sum_{y\in \La_c} f_\kappa (y/N) \right)^2}
{\mathrm{var}\left(\sum_{y\in \La^G_c} \gamma_N(y) f_\kappa (y/N) \s_y \right)}=
\frac
{(\int_D f_\kappa (x) \dd x)^2}
{\int_D \int_D f_\kappa (x)f_\kappa (x^\prime)R_d\vert x -
x^\prime\vert^{-d+2}\dd x\dd x^\prime}\equiv\mathrm{C}(f_\kappa),
\end{equation} 
where $R_d$ is defined in \eqref{eq:Rd}.
We stress once again 
that the variance appearing on the left--most term of
\eqref{eq:dil_cap} depends on $\sg$: but for $\sg\in E_N^\complement$
and fixed $f$
this convergence is uniform.
This tells us that for every $\eps>0$ and for every sufficiently large $N$
\begin{equation}
\sup_{\sg \in E_N^\complement}
\bPr\otimes \Prob
\left( F_N\big \vert \mathcal{F}^\sg\right)(\sg)
\le \exp \left(-N^{d-2} \log N \left((1-4\eta)^2 \alpha \mathrm{C}(f_\kappa)
+\eps\right)/2\right), \end{equation}
and recalling \eqref{eq:UBviaFandE} we obtain
\begin{equation}
\limsup_{N\rightarrow \infty}
\frac1{N^{d-2}\log N}
\log\bPr \otimes \Prob (\Omega_{N,\sigma}^+) \le 
-(1-4\eta)^2 \alpha \mathrm{C}(f_\kappa)/2.
\end{equation}  
We can then let $\alpha\nearrow 4(G_L+Q)$,
$L\nearrow \infty$, $\eta \searrow 0$ and 
$\kappa\searrow 0$ and optimise over the choice of    
$f_\kappa$, which is
now any function which is piecewise constant  over an arbitrarily thin regular
grid and equal to zero outside $D$. By the second line in 
\eqref{eq:capacity}, the capacity of $D$ appears and 
 we are done.

\qed

\bigskip

We complete this section 
by observing that 
one can extract from an annealed upper bound 
on the probability of $\Omega_{N,\sg}^+$ a quenched upper
bound. The annealed upper bound is provided 
in Proposition~\ref{th:UBsgm}.

\bigskip
\begin{cor}
\label{th:UBquenched}
Assume H.\ref{hp:hyp1}, 
 H.\ref{hp:hyp2} and H.\ref{hp:hyp3},
we have that
\begin{equation}
\label{eq:UBquenched}
\limsup_{N\rightarrow \infty}
\frac1{N^{d-2}\log N}\log \Prob
\left( \Omega_{N,\sigma}^+ \right)\le -2(G+Q) \mathrm{Cap}(D), 
\ \ \ \ \ \ \ \ \ \   \bPr (\dd \sg)\text{--a.s.}.
\end{equation}
\end{cor}

\bigskip
\noindent
{\it Proof.}
Let $X_N(\sg)$ be the random variable
$\log \Prob (\Omega_{N,\sigma}^+ )/(N^{d-2}\log N)$
and choose $\ell =2(G+Q) \mathrm{Cap}(D)-\eps$, $\eps>0$.
By Markov inequality  we have that
\begin{equation}
\label{eq:fromMarkov}
\frac1{N^{d-2}\log N}
\log\bPr \left(X_N \ge -\ell \right)
\le \ell  +  \frac1{N^{d-2}\log N} \log
\bPr \otimes \Prob \left( \Omega_{N,\sg}^+\right).
\end{equation}
Taking now the $\limsup$ on both sides, 
we get that for sufficiently large $N$
\begin{equation}
\bPr (X_N\ge -\ell)\le \exp(-\eps N^{d-2}\log N).
\end{equation}
Thus by Borel-Cantelli I, for all $\eps>0$,
\begin{equation}
\bPr(X_N \ge -2(G+Q)\text{Cap}(D)+\eps
\hbox{ i.o.})=0,
\end{equation}
whence the thesis.
\qed

\bigskip

\section{Entropic repulsion} 
\label{sec:path} 
\setcounter{equation}{0} 

This section is devoted to the proof 
of Theorem~\ref{th:height}. It is {\sl roughly}
split into two parts (lower and upper bounds)
even if some of the arguments require both
upper and lower estimates at the same time. 

\medskip
\subsection{Lower bounds}
We need the following
 preliminary result on the hitting probabilities
of the  simple random walk on $\Z ^d$. 
We denote by $\{ X_j^x\}_{j\ge0}$
 the simple random walk for which $X_0=x$, and its law by $P^x$.

\bigskip 

\begin{lem}\label{th:hitting} For any positive integer $n$ 
let $S_n=\{y\in\Z^d: \, \vert y\vert \ge n$ and  there exists
$x\in \Z^d$ such that $\vert x\vert <n$ and $\vert y-x\vert=1\}$,
 $\tau_x=\inf \{ j\ge 0: 
X_j^x \in S_n\}$ and 
\begin{equation}
H (x,y)=P^x \left( X_{\tau_x}=y\right),
\end{equation}
for $\vert x\vert<n$ and $y\in S_n$.
Then there exist $c_1,c_2,c_3>0$ such that 
for every  $\epsilon \in (0,1/4)$
\begin{equation}
\begin{split}
& c_1 n^{1-d}\le H(x,y)\le c_2 n^{1-d},\cr
& \vert H(x,y)-H(x^\prime ,y)\vert \le c_3\eps n^{1-d},
\end{split}
\label{eq:hitting}
\end{equation}
for every  $x,x^\prime$ such that $\Vert x \Vert \vee \Vert x^\prime \Vert \le \epsilon n$   
and every $y\in S_n$.
\end{lem}

\bigskip
\noindent
{\it Proof.} In   \cite[Lemma 1.7.4]{cf:Lawler}
it is shown that
\begin{equation}
\label{eq:ublbonhp}
c_1 n^{1-d}\le H(0,y)\le c_2 n^{1-d},
\end{equation}
and from the proof of Theorem 1.7.1 in 
\cite{cf:Lawler}, where the author proves that
for every fixed $u$ one can find
a positive constant $c$ such that   
$
\vert H(u,y)-H(0,y)\vert \le c_u O(n^{-1})n^{1-d}
$,
it is not difficult to see that 
one can choose $c_u=c_3\vert u\vert$, for some fixed constant $c_3$, so that
if $\Vert x \Vert \le \epsilon n$
\begin{equation}
\label{eq:ubondiffhp}
\vert H(x,y)-H(0 ,y)\vert \le c_3 \epsilon n^{1-d}.
\end{equation}
By combining 
\eqref{eq:ublbonhp} and 
\eqref{eq:ubondiffhp}, possibly redefining $c_1$, $c_2$ and $c_3$,
we get \eqref{eq:hitting}.
\qed

\bigskip

For what follows it turns out to be 
convenient to introduce the notion of {\sl empirical measure}:
given $A$ finite subset of $\Z^d$ and $I\subset \R$
we define that function $L_A (I): \R^{\Z^d} \to [0,1]$
as 
\begin{equation}
(L_A (I))(\s)=\frac1{\vert A\vert} \sum_{x\in A} \ident_{I} (\s_x).
\end{equation} 
If $I$ is an interval, say $I=(a,b\,]$, then we 
drop the extra parentheses: $L_A(I)=L_A(a,b\,]$.
The main result of this subsection is the following:

\bigskip

\begin{pro}
\label{th:lbhempav}
For any $a<4(G+Q)$ and every $\delta>0$
\begin{equation}
\label{eq:lbhempav}
\lim_{N\to \infty}
\Prob \left( L_{D_N}(-\infty, \sqrt{a\log N})\ge \delta \big\vert
\Omega^+_{N,\sg}\right)=0, 
\ \ \ \ \ \bPr (\dd  \sg)\text{--a.s..}
\end{equation}
\end{pro}

\bigskip

\noindent
{\it Proof.}
We 
adopt the notation
of Section \ref{sec:UBrsannealed}. The essential difference
here is that $L$ is not a fixed (large) number: rather we 
choose $L=L(N)\nearrow \infty$ as $N \nearrow \infty$.
In what follows $\epsilon $ is a small positive number,
that we will choose in the last steps of the proof,
and we use the short--cut notation 
$B_{\epsilon L} (y)=B_{\lfloor\epsilon L \rfloor}(y)$ and
$D_N^\epsilon= \bigcup _{y\in \Lambda_c} B_{\epsilon L} (y)$.
Of course $\lim_{N\to \infty}
\vert D_N^\epsilon \vert/ \vert D_N\vert =\epsilon^{d}$.

We start with the following remark: it suffices to prove
that for every $\delta>0$ there exists $\epsilon>0$ such that
\begin{equation}
\label{eq:lbhempav1}
\lim_{N\to \infty}
\Prob \left( L_{D^\epsilon_N}(-\infty, \sqrt{a\log N})\ge \delta \big\vert
\Omega^+_{N,\sg}\right)=0, 
\ \ \ \ \ \bPr (\dd  \sg)\text{--a.s..}
\end{equation}
In fact the full result, i.e. \eqref{eq:lbhempav},
is a direct consequence of a finite number (approximately $\epsilon^{-d}$)
of repetitions of the same argument applied to shifted copies of $\Lambda_c$.

We prove two lemmas with which we select
$L=L(N)$ and a {\sl good subset} of $\Lambda_c$: note that
these two lemmas concern $\bPr$ and not $\Prob$.

\bigskip

\begin{lem}
\label{th:lem1sg}
For every $\varrho \in (0,2)$ and  $\zeta <2Q\varrho$
choose $L=2\lfloor N^{\varrho/d}\rfloor$. Then
$\bPr(\dd \sg)$--a.s.~there exists $N_0(\sg)<\infty$
such that
for every  $N>N_0(\sg)$ the following  holds:
for every $y\in \Lambda_c$ there exists $\widetilde x (y)\in
B_{\epsilon L} (y)$ such that $\sg_{\widetilde x (y)}
\ge \sqrt{\zeta \log N}$.
\end{lem}

\bigskip
\noindent
{\it Proof.} 
Set $E_N=\{\sg:$ there exists $y\in \Lambda_c$ such that
$\sg_x < \sqrt{\zeta \log N}$ for every $x \in B_{\epsilon L} (y)\}$.
We need to show that $\bPr (E_N\ \text{i.o.})=0$.
Since the $\sg$--field is IID we have that for sufficiently large $N$ 
\begin{equation}
\begin{split}
\bPr (E_N)&=
1-\left(
1-\bPr \left(
\sg_x <\sqrt{\zeta \log N} , \text{ for every } x \in B_{\epsilon L}(0)
\right)
\right)^{\vert \Lambda_c\vert}\\
&\le 1-\left(1-p_N^{c_1N^\varrho}\right)^{c_2 N^{d-\varrho}},
\end{split}
\end{equation}
where $c_1$ and $c_2$ are positive constants
and $p_N= 1-N^{-(\zeta/2Q)-\eps^\prime}$:
$\eps^\prime$ is any strictly positive real number (we have applied H.\ref{hp:hyp2}) 
that we choose smaller than $(\varrho/2)-
\zeta/(4Q)$. We conclude that  
$\bPr (E_N)\le \exp(-N^{(\varrho/2)-
\zeta/(4Q)})$ for sufficiently large $N$ and therefore,
by Borel--Cantelli I, the proof is complete.
\qed (Lemma~\ref{th:lem1sg}).

\bigskip
For the second lemma we need some notation: set $S_L(y)=
\{z: \vert z-y \vert \ge L$ and there exists $x$ such that
$\vert x-y\vert<L$ and $\vert x-z\vert=1\}$. For any $\s \in \R^{\Z^d}$
and every $x \in B_{\epsilon L}(y)$
define $M_x^\circ (\s) = \sum_{z \in S_L(y)} H(x,z) \s _z$
(note that this $M_x^\circ(\s)$ is different from $M_x^\Box(\s)$ as defined in
the previous section).

\bigskip

\begin{lem}
\label{th:lem2sg}
Let $a>0$ and choose $L=L(N)$ such that $\lim_{N\to \infty} L(N)/N^q $ is positive
and finite for a given  $q\in (0,1)$. Then for every $\delta>0$ there exists
$\epsilon_0>0$ such that 
$\bPr (\dd \sg)$--a.s.~there 
exist $\widetilde{\Lambda_c^G}\subset \Lambda_c$ and a finite number $N_0$
satisfying the following properties:
\begin{itemize} 
\item 
for every choice of $y \in \widetilde{\Lambda_c^G}$ and every $\epsilon\in
(0,\eps_0)$, if $\s \in \Omega^+_{N,\sg}$ and if there exists $x \in
B_{\epsilon L} (y)$ such that $M_x^\circ (\s) \le \sqrt{a\log N}$ then
\begin{equation}
\max_{x^\prime \in B_{\epsilon L}(y)} \left \vert M_{x^\prime}^\circ (\s)
-  M_{x}^\circ (\s) \right\vert
\le \delta \sqrt{ \log N}, 
\end{equation}
for every $N\ge N_0$,
\item
\begin{equation}
\lim_{N\to \infty} \frac{\vert \widetilde{\Lambda_c^G} \vert} {\vert
\Lambda_c\vert}=1. \end{equation}
\end{itemize}
\end{lem}

\bigskip
\noindent
{\it Proof.}
Let $y \in \Lambda_c$.
By applying repeatedly Lemma~\ref{th:hitting} we obtain
\begin{equation}
\begin{split}
\left \vert M_{x^\prime}^\circ (\s)
-  M_{x}^\circ (\s) \right\vert 
&
\le \sum_{z\in S_L(y)}
\left \vert H(x^\prime,z)-H(x,z)\right\vert
\left\vert \s_z \right\vert
\le c_3c_1^{-1}\epsilon \sum_{z\in S_L(y)}
H(x,z)\vert \s_z \vert
\\
& =c_3c_1^{-1}\epsilon M_x^\circ(\s)
 +2c_3c_1^{-1}\epsilon \sum_{z\in S_L(y):\s_z < 0}
H(x,z)\left\vert \s_z \right\vert
\\
& 
\le 
 c_3c_1^{-1}\epsilon \sqrt{a \log N}+
\frac{2c_3^2c_1^{-1}\epsilon c}
{\left\vert S_L(y)\right\vert}
\sum_{z\in S_L(y)} \left \vert \sg_z \right\vert.
\end{split}
\end{equation}
The Lemma is therefore proven once we show
for example that there exists a sequence $\{\delta_N\}_N$ of
positive numbers, $\delta_N\searrow 0$ as $N\nearrow \infty$,
such that
\begin{equation}
\label{eq:BC1sglbh}
\sum _N \bPr
\left(
\left\vert 
  \left\{ 
      y\in \Lambda_c:
      \frac1{\left\vert S_L(y)\right\vert}
      \sum_{z\in S_L(y)} 
       \left\vert \sg_z \right\vert >\frac{\sqrt{a\log N}}{2c_3}
  \right\}
\right\vert
> \delta_N \left\vert \Lambda_c \right\vert
\right)
<\infty.
\end{equation}
Choose 
\begin{equation}
\delta_N^2=\max
\left(
\bPr 
\left(
           \frac1{\left\vert S_L(y)\right\vert}
            \sum_{z\in S_L(y)} \left \vert \sg_z \right\vert >\frac{\sqrt{a\log N}}{2c_3}
\right), \left\vert \Lambda_c \right\vert^{-1/2}
\right).
\end{equation}
By H.\ref{hp:hyp3} $\delta_N$ vanishes as $N$ tends to infinity.
The rest of the proof of \eqref{eq:BC1sglbh}
follows from a direct application of \eqref{eq:Bennett}.
 \qed (Lemma~\ref{th:lem2sg})

\bigskip

We now choose $\sg$ in the {\sl good set}
specified by Lemma~\ref{th:lem1sg} and Lemma~\ref{th:lem1sg}.
Let us fix the choice of the parameters with the help
of an extra  parameter $\widetilde{\eps}>0$:
\begin{equation}
\label{eq:manyparam}
\varrho =(2Q/(G+Q))-\widetilde{\eps}>0, \ \ \ 
\zeta =
\left(\sqrt{\frac{2Q}{\sqrt{G+Q}}}-\widetilde{\eps}\right)^2>0, 
\ \ \
L=2\lfloor N^{\varrho/d}\rfloor.
\end{equation}

Going back to the proof of Proposition~\ref{th:lbhempav} let us make
another step in the spirit of \eqref{eq:lbhempav1}: 
we claim in fact that \eqref{eq:lbhempav1}  is proven 
if we show that there exists $\widetilde{\delta}>0$ such that
for every $x\in B_{\epsilon L}(0)$ 
\begin{equation}
\label{eq:lbhempav2}
\lim_{N\to \infty}\frac{1}{N^{\widetilde{\delta}}}\log
\Prob \left( L_{\widetilde{\Lambda^G_c} +x}(-\infty, \sqrt{a\log N})\ge \delta
\big\vert
 \Omega^+_{N,\sg}\right)<0, 
\ \ \ \ \ \bPr (\dd  \sg)\text{--a.s..}
\end{equation}
The claim follows since 
\begin{equation}
\left\{ \s: L_{D^\epsilon_N}(-\infty, \sqrt{a\log N})\ge \delta\right\}
\subset
\bigcup_{x\in B_{\epsilon L}}
\left\{ \s :
L_{\Lambda_c +x}(-\infty, \sqrt{a\log N})\ge \delta \right\},
\end{equation}
and $\vert B _{\epsilon L}\vert \exp(-cN^{\widetilde \delta})$ vanishes
as $N \to \infty$.
Lemma~\ref{th:lem2sg} guarantees that we may substitute $\Lambda_c$ with
$\widetilde{\Lambda_c^G}$.

We think now of $x$ as fixed and observe that $M_y^\circ(\s)$ and $\s_y$, for
$y\in \widetilde{\La^G_c} +x$ are close in the sense specified by the following
lemma.

\bigskip

\begin{lem}
\label{th:harmmeansands}
For every $\eta>0$, $\delta^\prime>0$,
\begin{equation}\label{eq:hmean1}
\limsup_{N\to\infty} \frac 1{N^{d-\varrho}}
\log\Prob \left(
\left \vert \left\{y\in  \widetilde{\La^G_c} +x\,:\,
\vert\s_y-M_y^\circ (\s) \vert>\eta\sqrt{\log N}\right\} \right \vert 
\ge \delta^\prime\vert \widetilde{\La^G_c} \vert
\right)<-c,
\end{equation} 
for some $c>0$.
\end{lem}
\bigskip

\noindent
{\it Proof.} We observe that
 $\{ \s_y-M_y^\circ (\s) \}_{y\in \widetilde{\Lambda_c^G}+x}$ forms an
IID collection of centered Gaussian random variables of variance that is
not larger than $G$. Therefore for every $\eta>0$
\begin{equation}
\limsup_{N\to \infty}
\frac1{N^{d-\varrho}}
\log 
\Prob \left(
\frac{1}{\vert \widetilde{\Lambda_c^G}\vert}
\sum_{y\in \widetilde{\Lambda_c^G} +x } \left \vert \s_y-M_y^\circ (\s) 
\right \vert  > \eta \sqrt{\log N}
\right)<0.
\end{equation} 
Now note that the probability in \ref{eq:hmean1} is
not larger than
\begin{equation}
\Prob \left(
\frac{1}{\vert \widetilde{\Lambda_c^G}\vert}
\sum_{y\in \widetilde{\Lambda_c^G} +x } \left \vert   \s_y-M_y^\circ (\s) 
\right \vert  \ge \eta\delta^\prime \sqrt{\log N}
\right),
\end{equation}
whence the thesis. \qed (Lemma~\ref{th:harmmeansands})

\bigskip
We define
\begin{equation}
E_N^{\widetilde\eps}=\left\{
\left \vert \left\{y\in  \widetilde{\La^G_c} +x\,:\,
\vert\s_y- M_y^\circ(\s)\vert \le\widetilde\eps\sqrt{\log N}\right\} \right
\vert  > (1- \delta/2)\vert \widetilde{\La^G_c} \vert
\right\}.
\end{equation}
By Lemma~\ref{th:harmmeansands} we know that $\Prob({E_N^{\widetilde\eps}}^\complement)
<\exp(-cN^{d-\varrho})$ for sufficiently large $N$ and for some positive $c$.
In order to prove \eqref{eq:lbhempav2}, 
we analyze
\begin{equation}
\label{eq:inters1}
\begin{split}
\Prob & \left(
\left \vert \left\{y\in  \widetilde{\La^G_c} +x\,:\,
\s_y<\sqrt{a\log N}\right\} \right \vert 
\ge \delta\vert \widetilde{\La^G_c} \vert;\, \Omega_{N,\sigma}^+
\right)\\
& = \Prob \left(
\left \vert \left\{y\in  \widetilde{\La^G_c} +x\,:\,
\s_y < \sqrt{a\log N}\right\} \right \vert 
\ge \delta\vert \widetilde{\La^G_c} \vert;\,{E_N^{\widetilde\eps}}^\complement;
\, \Omega_{N,\sigma}^+ \right)\\
&+
\Prob \left(
\left \vert \left\{y\in  \widetilde{\La^G_c} +x\,:\,
\s_y <\sqrt{a\log N}\right\} \right \vert 
\ge \delta\vert \widetilde{\La^G_c} \vert; \, E_N^{\widetilde\eps}; \,
\Omega_{N,\sigma}^+ \right)
\end{split}
\end{equation}
The first term in the right--hand side of \ref{eq:inters1} is
not larger than
 $\exp(-cN^{d-\varrho})$.
We focus on the second term:
\begin{equation}\label{eq:rhs1}
\begin{split}
\Prob & \left(
\left \vert \left\{y\in  \widetilde{\La^G_c} +x\,:\,
\s_y < \sqrt{a\log N}\right\} \right \vert 
\ge \delta\vert \widetilde{\La^G_c} \vert; \,
E^{\widetilde\eps}_N;\,\Omega_{N,\sigma}^+ \right)\\
& \le
\Prob \left(
\left \vert \left\{y\in  \widetilde{\La^G_c} +x\,:\,
M_y^\circ(\s) < \sqrt{a\log N}+\widetilde\eps\sqrt{\log N}
\right\} \right \vert 
\ge \delta\vert \widetilde{\La^G_c} \vert/2; \, \Omega_{N,\sigma}^+
\right).
\end{split}
\end{equation}
Now we use the fact  that,
 by Lemma~\ref{th:lem2sg}, when $M_y^\circ(\s)\le\sqrt{b\log N}$ for
some $b>0$,
$\vert M_{\widetilde x(y)}^\circ(\s)-
M_y^\circ(\s)\vert \le \widetilde\eps\sqrt{\log N}$ if we choose  $\eps$
sufficiently small  and $N$ sufficiently
large. We recall that $\eps$ was introduced at the beginning of the
proof.
Hence the last term in \ref{eq:rhs1} is not larger than
\begin{equation}\label{eq:rhs2}
\Prob \left(
\left \vert \left\{y\in  \widetilde{\La^G_c} +x\,:\,
M_{\widetilde x(y)}^\circ(\s) < (\sqrt{a}+2\widetilde\eps)\sqrt{\log N}
\right\} \right \vert 
\ge \delta\vert \widetilde{\La^G_c} \vert/2; \, \Omega_{N,\sigma}^+
\right).
\end{equation}
Set 
\begin{equation}
\widetilde{E}_N=\left\{
\s:\left \vert \left\{y\in  \widetilde{\La^G_c} +x\,:\,
M_{\widetilde x(y)}^\circ(\s) < (\sqrt{a}+2\widetilde\eps)\sqrt{\log N}
\right\} \right \vert 
\ge \delta\vert \widetilde{\La^G_c} \vert/2
\right\},
\end{equation}
with this notation the last term in \eqref{eq:inters1} is
dominated by
\begin{equation}
\label{eq:pomega1}
\begin{split}
\Prob(\widetilde{E}_N\cap \Omega^+_{N,\sigma}) & \le \E\left[ \prod_{y\in
\widetilde{\La^G_c}}
\Prob\left( \s_{\widetilde x(y))} \ge\sigma_{\widetilde x(y)}
\vert\cF^\s_{S_L(y)} \right); \, \widetilde{E}_N\right]\\
&\le \E\left[ \prod_{y\in\widetilde{\La^G_c}}
\Prob\left( \s_{\widetilde x(y)} \ge\sqrt{\zeta\log N} \vert\cF^\s_{S_L(y)}
\right); \, \widetilde{E}_N\right].
\end{split}
\end{equation}
But on $\widetilde{E}_N$ we have that (we may think $a\ge 2Q/\sqrt{G+Q}$)
\begin{equation}
\begin{split}
 \prod_{y\in\widetilde{\La^G_c}}
\Prob\left( \s_{\widetilde x(y)} \ge\sqrt{\zeta\log N} \vert\cF^\s_{S_L(y)}
\right)
&\le
\left(
1-\Phi\left(-\frac{(\sqrt a+2\widetilde \eps -\sqrt\zeta)
\sqrt{\log N}}{\sqrt{G_L}}
\right)\right)^{\delta\vert\widetilde{\La_c^G}\vert/2} \\
&\le \exp\left(-N^{d-\varrho-\frac{(\sqrt a+3\widetilde \eps-\sqrt\zeta)^2}
{2G}}\right),
\end{split}
\end{equation}
where 
$G_L=\min\left\{{\mathrm var}\left(\s_x\vert\cF^\s_{S_L(y)}\right):
x\in B_{\eps L}(y)\right\}$ (note that $G_L$ does not depend on $y\in\La_c$
and $G_L\nearrow G$ as $L\nearrow\infty$). The last step
holds for sufficiently large $N$.

Recalling the estimate of $\Prob(\Omega^+_{N,\sigma})$
(Theorem~\ref{th:main1}),
equation~\ref{eq:lbhempav2} follows if
\begin{equation}
\label{eq:parineq}
d-\varrho-\frac{(\sqrt a+3\widetilde\eps-\sqrt\zeta)^2}
{2G}>d-2.
\end{equation}
A straightforward computation shows that if $\widetilde \epsilon<(\sqrt{4(G+Q)}-a)/4)$
then the left--hand side of \eqref{eq:parineq} is bounded
below by $d-2+\widetilde\epsilon$ and we are done.

\qed

\bigskip
\subsection{Upper bounds}
For $\Lambda \subset \R^d$,  $N\in \N$ and $\s\in \R^{\Z^d}$ 
we set $M_N^\La (\s)=\sum_{x\in \La_N}
\s _x / \vert \La_N \vert$: we 
 always consider
$\Lambda$ a bounded open set with piecewise 
smooth boundary (even if this condition could
be very much relaxed).
We observe that
$M_N(\s) \sim {\mathcal N}(0, N^{2-d}(c(\La)+o(1)) )$,
where $c(\La)>0$.

We give the following upper bound on the path of the 
interface above the rough wall:

\bigskip
\begin{pro}
\label{th:UBlem}
For every $\La \subset D$ we have that
\begin{equation}
\label{eq:UBlem}
\limsup_{N\to \infty}
\frac{\E \left[M_N^\La (\s)\big\vert \Omega_{N,\sigma}^+
\right]}{\sqrt{\log N}}\le \sqrt{4(G+Q)},
\end{equation}
$\mathbf{P}(\dd \sigma)$--a.s..
\end{pro}
\bigskip

\noindent
{\it Proof.} Set $\Prob_N\equiv \Prob T_{\psi_N}^{-1}$,
${\psi_N}_x \in [0,\infty)$ independent of $x$.
We observe that  $\Prob(\cdot \vert\Omega_{N,\sigma}^+)$
is dominated by $\Prob_N
(\cdot \vert\Omega_{N,\sigma}^+)$.
This follows for example by writing a finite volume 
approximation of $\Prob$, with $0$--boundary conditions,
namely $\mu_n (\cdot)=\Prob (\cdot \vert \mathcal{F}_{{D_n}^\complement})
 (\psi^0)$,
$\psi^0\equiv 0$: we view this measure as a measure
on $\R^{D_n}$. One verifies directly 
that if $T: \R^{D_n}\to \R^{D_n}$ is defined by $(T \s)_x=\s_x+a$, $a\ge 0$, 
then
 $\mu_n T^{-1}$ dominates $\mu_n$ in the strong FKG sense
  (that is the two measures
satisfy Holley's inequality, cf. \cite{cf:Preston}). 
Therefore, if $n\ge N$, 
we can define
 $\mu_n T^{-1}(\dd \s)
\exp(-\sum_x U_x(\s_x))/Z$ and  $\mu_n (\dd \s)
 \exp(-\sum_x U_x(\s_x))/Z^\prime$, with $U_x (\cdot)$ a potential
that for definitness we choose equal to 
$\beta r^4 \ident_{(-\infty, \sigma_x]}(r) \ident_{D_N}(x)$, $\beta>0$ and 
$Z, Z^\prime$ are the normalization
constants,  and this two new measures
are still ordered in the strong FKG sense. The limit for $n \to \infty$
and then 
$\beta \to \infty$ 
recovers the desired inequality.

We choose ${\psi_N}_x=\alpha_N $, with 
$\alpha_N=  \sqrt{(4(G+Q)+\epsilon)\log N}$ and 
$\epsilon>0$.
Therefore
\begin{equation}
\begin{split}
\label{eq:FKGandcomp}
\E \left[ M_N^\La (\s)\big\vert \Omega_{N,\sigma}^+
\right]
&\le \E_N \left[ M_N^\La (\s)\big\vert \Omega_{N,\sigma}^+
\right]
\\
&=\alpha_N +
 \E \left[ M_N^\La (\s)\big\vert \Omega_{N,\sigma-\psi_N}^+
\right].
\end{split}
\end{equation}
By applying \eqref{eq:entineq2}
with $Y=\pm  M_N^\La (\s)$, $t=\delta N^{d-2}$ ($\delta>0$)
and $E= \Omega_{N,\sigma-\psi_N}^+$ 
we have that 
\begin{equation}
\left\vert
\E \left[ M_N^\La (\s)\big\vert \Omega_{N,\sigma-\psi_N}^+
\right]
\right\vert
\le 
\delta (c(\La)+o(1))-\frac1{\delta N^{d-2}}
\log
\Prob \left(\Omega_{N,\sigma-\psi_N}^+\right).
\end{equation}
In Section \ref{sec:LBrs}, Lemma~\ref{th:technheight}, we have shown 
that for
every $\epsilon>0$ 
\begin{equation}
\lim_{N\to \infty}
\frac1{N^{d-2}}
\log
\Prob \left(\Omega_{N,\sigma-\psi_N}^+\right)=0.
\end{equation}
$\mathbf{P}(\dd \sigma)$--a.s..
Since $\delta>0$ can be chosen arbitrarily small
we  conclude that 
\begin{equation}
\lim_{N\to \infty}
\E \left[ M_N^\La (\s)\big\vert \Omega_{N,\sigma-\psi_N}^+
\right]=0.
\end{equation}
This is more than we need: 
apply it in  \eqref{eq:FKGandcomp} to  get \eqref{eq:UBlem}.
\qed

\bigskip


\begin{pro}
\label{th:LBlem}
For every $\La \subset D$ we have that
\begin{equation}
\label{eq:LBlem}
\liminf_{N\to \infty}
\frac{\E \left[M_N^\La (\s)\big\vert \Omega_{N,\sigma}^+
\right]}{\sqrt{\log N}}\ge \sqrt{4(G+Q)},
\end{equation}
$\mathbf{P}(\dd \sigma)$--a.s..
\end{pro}
\bigskip

\noindent
{\it Proof.} Since Proposition~\ref{th:lbhempav} may

be proved with $\La$ in place of $D$,
for all positive $\eps$ and sufficiently large $N$,

there exists
$\Omega_\eps\sbs \R^{\Z^d}$ such that $\Prob \left(
\Omega_\eps\,\vert\,\Omega^+_{N,\sigma}\right)\ge (1-\eps)$,
for all $\s\in\Omega_\eps$ there exists $A_\eps\sbs\La_N$,
$\vert A_\eps\vert >(1-\eps)\vert \La_N\vert$ and 
$\s_x\ge\sqrt{(4(G+Q)-\eps)\log N}$ for every $x\in A_\eps$.

If $\s\in\Omega_\eps$, we decompose 
\begin{equation}\label{eq:decempmean}
M_N^\La(\s)=\frac 1 {|\La_N|}\sum_{x\in A_\eps}
	\s_x+\frac 1 {|\La_N|}\sum_{x\in \La_N\setminus A_\eps}\s_x.
\end{equation}
Under $\Prob(\cdot\vert\,\Omega_{N,\sigma}^+)$, the
first term in equation \ref{eq:decempmean}
is larger than $(1-\eps)\sqrt{(4(G+Q)-\eps)
\log N}$,
while the second term is larger than 
$- \sum_{x\in\La_N}|\sigma_x|/\vert\La_N\vert$.
Observe that the last quantity is also a minorant
for $M_N^\La(\s)$ when $\s\in\Omega_\eps^\complement$.
But $- \sum_{x\in\La_N}|\sigma_x|/\vert\La_N\vert$
converges ${\bf P}(\dd\sigma)$--a.s.~to 
${\bf E}[\vert\sigma_0\vert ]$, thus

\begin{equation}
\liminf_{N\to \infty}
\frac{
\E \left[M_N^\La (\s)\big\vert \Omega_{N,\sigma}^+
\right]}{\sqrt{\log N}}\ge 
(1-\eps)^2\sqrt{(4(G+Q)-\eps)
}.
\end{equation}
The thesis follows taking $\eps\searrow 0$.
\qed

\bigskip 
\begin{rem} 
\label{th:hydrolim}
We have therefore that, for every choice of $\Lambda$,
$\bPr (\dd \sg)$--a.s.
\begin{equation}
\label{eq:MLambdalimit}
\lim_{N\to \infty}
\frac{\E \left[M_N^\La (\s)\big\vert \Omega_{N,\sigma}^+
\right]}{\sqrt{\log N}}= \sqrt{4(G+Q)}.
\end{equation}
By the Brascamp--Lieb inequality \cite{cf:BL}  the random variable $M_N^\La (\s)$,
(even) under the conditioned measure $\Prob (\dd \s\vert \Omega_{N,\sigma}^+)$
\cite{cf:DG}, 
has a {\sl sub--Gaussian} behavior: the exponential centered moment
of $M_N^\La (\s)$ is bounded by $\exp(cN^{2-d})$.
This immediately yields the {\sl hydrostatic limit} of the field:
if for $r\in D$ we define $u_N(r)=\s_{rN}/\sqrt{\log N}$ for $rN\in \Z^d$
and if we extend $u_N$ to a function from $D$ to $\R$
(for example) by a polylinear interpolation, then
$u_N$ converges weakly to $u\equiv \sqrt{4(G+Q)}$ (that is $\int_D u_N f \to
\sqrt{4(G+Q)}\int f$ for every $f\in C^0_b(D;\R)$), 
in probability, with respect to $\Prob (\dd \s \vert \Omega_{N,\sigma}^+)$,
and $\bPr (\dd \sg)$--a.s.. 
\end{rem}

\bigskip

{\begin{pro}
\label{th:ubheight}
For any $b>4(G+Q)$ and every $\delta>0$, 
\begin{equation}
\lim_{N\to\infty}\Prob\left(L_{D_N} (\sqrt{b\log N},+\infty)
\ge\delta \vert\Omega_{N,\sigma}^+\right)\longrightarrow 0,
\end{equation}
${\bf P}(\dd\sigma)$--a.s..
\end{pro}

\bigskip
\noindent
\it Proof.}
Fix $b$ and define 
\begin{equation}
\overline{N_b}(\s)=\left\vert \{x\in D_N : \s_x>\sqrt{b \log N}\}\right\vert.
\end{equation} 
By Proposition~\ref{th:lbhempav}, for all positive $\eps$ and sufficiently 
large $N$,
there exists $\Omega_\eps$ such that
$\Prob\left(\Omega_\eps\vert\Omega^+_{N,\sigma}\right)>(1-\eps)$ and  on
$\Omega_\eps$, $\s_x\ge \sqrt{(4(G+Q)-\eps)\log N}$ on at
least $(1-\eps)\vert D_N\vert$ sites $x$. 
Thus on $\Omega_\eps$, $\s_x/\sqrt{\log N}$ is larger than $\sqrt b$ on at
least $\overline{N_b}(\s)$ sites, larger than $\sqrt{4(G+Q)-\eps}$ on
at least $(1-\eps)\vert D_N\vert -\overline{N_b}(\s)$ sites and  
on the remaining (at most
$\eps \vert D_N\vert$) sites it is larger than
$-\vert\sigma_x\vert/\sqrt{\log N}$ (thanks to the conditioning on
$\Omega^+_{N,\sigma}$). 
Thus
\begin{equation}
\frac{\E\left[
M_N^{D}(\s)\vert\Omega_{N\sigma}^+\right]}{\sqrt{\log N}}
>(1-\eps)f(b,\eps)-(1+\eps)\frac{\sum_{x\in D_N}\vert\sigma_x\vert}
{\sqrt{\log N}\vert D_N\vert},
\end{equation}
where 
\begin{equation}
f(b,\eps)= 
(1-\eps)\sqrt{4(G+Q)-\eps}+\E\left[\left.\frac{\overline{N_b}(\s)}
{\vert D_N\vert}
\right\vert\Omega_{N\sigma}^+\right]
(\sqrt b-
\sqrt{4(G+Q)-\eps}).
\end{equation}
Now we let $N$ grow to infinity: by Remark~\ref{th:hydrolim} we obtain
\begin{equation}
\sqrt{4(G+Q)} \ge (1-\eps)\sqrt{4(G+Q)-\eps}+
\limsup_{N\to\infty}\E\left[\left.\frac{
\overline{N_b}(\s)}{\vert D_{N} \vert}
\right\vert\Omega_{N\sigma}^+\right]
(\sqrt b-
\sqrt{4(G+Q)-\eps}).
\end{equation}
Now let $\eps\searrow 0$: since $b$ is chosen
strictly larger than $\sqrt{4(G+Q)}$ we get that 
\begin{equation}
\limsup_{N\to\infty}\E\left[\left.\frac{
\overline{N_b}(\s)}{\vert D_{N} \vert}
\right\vert\Omega_{N\sigma}^+\right]
=0.
\end{equation}
This leads to the conclusion, once we observe that
$L_{D_N}(\sqrt{b\log N},+\infty)=\overline{N_b}(\s)/
\vert D_N\vert$. \qed
\bigskip

\section{Super--Gaussian and sub--Gaussian regimes}
\label{sec:mainprime}
\setcounter{equation}{0} 

The proof of Theorem~\ref{th:mainprime}
can be obtained following and modifying step by step 
the proof of Theorem\ref{th:main1} and  
Theorem\ref{th:height}. However 
many of the steps in such an approach would be superfluous:
we therefore sketch the proof pointing out the most
substantial simplifications. On the way we 
also give some results that sharpen Theorem~\ref{th:mainprime}.
We assume H.\ref{hp:hyp1} and H.\ref{hp:hyp3}.
\medskip

\subsection{The sub--Gaussian regime}
Under H.2--\ref{hp:hyp2-1} one immediately sees that
for every $\theta>0$ and every $k>0$
\begin{equation}
\label{eq:maxsuperG}
\lim_{N\to \infty}
N^k \mathbf{P} \left(\max_{x\in D_N}
\sg _x >\theta \sqrt{\log N} \right)=0.
\end{equation}
Therefore in proving the lower bound
corresponding to \eqref{eq:quenannprime1} we may substitute the 
auxiliary field $\widetilde{\sg}$, previously
defined in \eqref{eq:aux}, with
$\tilde{\sg}_x = \theta \sqrt{\log N}$, with $\theta$
arbitrarily small.
At this point we may directly apply \cite[Prop.~2.1]{cf:BDZ}, that
is the lower bound in the case of a flat wall: by sending $\theta $
to zero we obtain the result.

For what concerns a proof of the upper bound corresponding to 
 \eqref{eq:quenannprime1}, due to the
weakness of the assumption H.\ref{hp:hyp3} the results in
\cite{cf:BDZ} are no longer applicable and one need some argument
in the spirit of the proof of 
\eqref{th:UBsgm}. The guideline is the following:
leave Definition \eqref{eq:Eea}  unchanged
and replace \eqref{eq:twofc} with
\begin{equation}
\prod_{y \in  \Lambda_c}  
\bPr \otimes \Prob \left(\s_y \ge \sigma_y\big \vert \cF^{\sg, \s } _{B(y)}\right) 
\le
\left( 1-
N^{-\frac{\alpha+\eps}{2 G_L}
}\right)^{\eta\vert C_{r}\vert}.
\end{equation}
The rest of the steps are identical (set $Q=0$).

The very same arguments apply in extending the
proof of Theorem~\ref{th:height}
 to cover  
the second part of Theorem~\ref{th:mainprime}(1).\qed
\medskip

We  stress that
observation \eqref{eq:maxsuperG}, which  allows a natural 
comparison argument, 
together with the result \eqref{eq:BDZDG} immediately yields 
the following sharpening of Theorem~\ref{th:mainprime}(1):

\bigskip
\begin{pro}
Under Hypothesis H.2--\ref{hp:hyp2-1} we have that
\begin{equation}
\limsup_{N\to \infty}
\sup_{x\in D_N}
\frac{\E \left[ \s_x \big \vert \Omega^+_{N,\sg}\right]}
{\sqrt{4G \log N}}\le 1.
\end{equation}
Moreover if also $-\sg$ satisfies 
H.2--\ref{hp:hyp2-1} then
\begin{equation}
\limsup_{N\to \infty}
\sup_{x \in D_N}
\left \vert
\frac{\E \left[ \s_x \big \vert \Omega^+_{N,\sg}\right]}
{\sqrt{4G \log N}} -1
\right \vert =0.
\end{equation}  
\end{pro}
\bigskip

\medskip

\subsection{The super--Gaussian regime}
Once again a look to the proof of \eqref{eq:LBsgm}
is sufficient to understand that the {\sl multiscale}
decomposition in  \eqref{eq:aux}
is superfluous. In proving the lower bound of
\eqref{eq:quenannprime2}
 we may substitute \eqref{eq:aux} with 
the much rougher discretization
\begin{equation}
\widetilde{\sigma}_x=
\begin{cases}
( (4Q+\epsilon)\log N)^{1/(2\beta)} &
\text{ if } \sg _x \le ( (4Q+\epsilon)\log N)^{1/(2\beta)},
\\
( (2d+2)Q\log N)^{1/(2\beta)}
 &
\text{ if } \sg _x \in ( ( (4Q+\epsilon)\log N)^{1/(2\beta)},(
(2d+2)Q\log N)^{1/(2\beta)} ],
 \\
\infty &\text{ otherwise,}
\end{cases}
\end{equation}
with $\epsilon >0$.
The rest of the proof of the lower bound for  
\eqref{eq:quenannprime2} follows in an absolutely
analogous, but simpler, way as the proof of
\eqref{eq:LBsgm}: the optimization over the
levels of the $\widetilde{\sg}$--field is
trivial.

For what concerns the upper bound for \eqref{eq:quenannprime2}
it suffices to redifine $E_{\eta, \alpha}$, cf. \eqref{eq:Eea}
in the proof of Proposition~\ref{th:UBsgm},
in the following way:
\begin{equation} 
\label{eq:Eeaprime}
E_{\eta, \alpha}= 
\left\{(\sg, \s): \text{there exists }r \in I \text{ such that } 
\vert \{ y \in C_r: M_y^\Box (\s )\le (\alpha\log N)^{1/2\beta}\}\vert 
\ge \eta \vert C_r \vert  \right\}. 
\end{equation} 
and one obtains (compare with \eqref{eq:twofc})
\begin{equation}
\label{eq:twofcprime} 
\prod_{y \in  \Lambda_c}  
\bPr \otimes \Prob \left(\s_y \ge \sigma_y\big \vert \cF^{\sg, \s } _{B(y)}\right) 
\le
\left( 1-
N^{-\frac{\alpha+\eps}{2Q}
}\right)^{\eta\vert C_{r}\vert}.
\end{equation}
The rest of the proof is essentially identical:
just substitute $\sqrt{\alpha\log N}$ with $(\alpha\log N)^{1/2\beta}$.
Analogous modifications to the proof
of Theorem~\ref{th:height}
completes the proof of Theorem~\ref{th:mainprime}(2).
\qed
\medskip

\begin{rem}
Since there 
are {\sl several} spikes
of the $\sg$--field going beyond the level of the interface,
in fact $\max_{x\in D_N} \sg_x \approx (2dQ\log N)^{1/2\beta}$,
one cannot hope to have a bound of the
type $\E [\s_x \vert \Omega_{N,\sg}^+] \le ((4Q+\epsilon)\log N)^{1/2\beta}$
uniformly in $x$ for $\epsilon $ arbitrarily small.
\end{rem}

\medskip
\begin{rem}
A word on heavier $\sg$--tails is due:
a new phenomenon is expected to arise if $\sg$ has 
{\sl power law} upward tails. In this case
$\max_{x\in D_N} \sg_x$, suitably normalized,
converges to a nondegenerate random variable and this
is sharply different of what happens in all the
 cases that we considered. Moreover excursion 
of the $\sg$--field beyond the level
$N^{\delta}$, some $\delta>0$, would now be possible,
even on more than $N^{\epsilon}$ sites, for some $\epsilon>0$
depending on $\delta$ and on the tail behavior.
It is quite clear from an entropy argument
that, even if $\epsilon <d-2$,
these spikes may have a very strong effect on the field:
almost local defomations of the $\s$--field 
are not
the optimal strategy to accomodate the presence
of the wall. This is in stark contrast with
the situation we dealt with, since (roughly) excursions of 
the $\sg$--field beyond level $(\log N)^k$
on $o(N^{d-2})$ sites produce only {\sl almost local}
modifications of the $\s$--field. 
\end{rem}

\bigskip

\section*{Acknowledgements}
We would like to thank Yvan~Velenik, who took part in the earlier
stages of this work, for many valuable 
exchanges.
G.G. would like to thank also Thierry~Bodineau and 
Ofer~Zeitouni for
important discussions on various
aspects of the proofs.
Part of this work has been developed while we were
 visiting the {\sl Institut Henri Poincar\'e} (fall 2001:
semester {\sl Limites Hydrodynamiques}):
we thank all the people at the institute for
the warm hospitality.


\bigskip


\begin{thebibliography}{15} 


\bibitem{cf:Abraham}
D.B.~Abraham, \textit{Surface structures and phase transitions--exact
results}, Phase transitions and critical phenomena, \textbf{10}, 1--74, Academic Press,
London, 1986.
 
\bibitem{cf:Bennett}
        G.~Bennett, \textit{Probability inequalities for the sum 
	of independent random variables}, 
	J. Amer. Stat. Assoc. \textbf{57}(1962), 33--45.

\bibitem{cf:Bolthausen}
E.~Bolthausen, 
\textit {Random walk representations and entropic repulsion for gradient
models}, 
Infinite Dimensional Stochastic Analysis, Koninklijke
Nederlandse Akademie van Wetenschappen, Ph. Cl\'ement et al. eds., (2000), 55--84.         

\bibitem{cf:BD}
E.~Bolthausen and  J.D.~Deuschel,
\textit{Critical large deviations for {Gaussian }
fields in the phase transition regime},
 Ann. Prob. \text\bf{21} (1994),
 1876--1920.  

\bibitem{cf:BDG} 
        E.~Bolthausen, J.-D.~Deuschel and G.~Giacomin, 
\textit{Entropic repulsion and 
the maximum of the two dimensional free field}, 
Ann. Probab. \textbf{29} (2001). 


\bibitem{cf:BDZ} 
        E.~Bolthausen, J.-D.~Deuschel and O.~Zeitouni, 
\textit{Entropic repulsion of the lattice free field}, 
Commun.~Math.~Phys. \textbf{170} (1995), 417--443. 

\bibitem{cf:BL}
 H.J.~Brascamp and E.~Lieb,
\textit{On extensions of the {B}run--{M}inkowski and
{P}rekopa--{L}einler theorems},
 J. Funct. Anal.
\textbf{22} (1976),
 366-389.

 \bibitem{cf:BeMF}
 J.~Bricmont, A.~el Mellouki and J.~Fr\"ohlich,
 \textit{Random surfaces in statistical mechanics:
 roughening, rounding, wetting},
  J. Stat. Phys. \textbf{42} (1986),
 743--798.



\bibitem{cf:DG} 
        J.-D.~Deuschel and G.~Giacomin, 
	\textit{Entropic repulsion
	for the free field: pathwise characterization in $d\ge3$},
	Commun.~Math.~Phys. \textbf{206}, 447--462 	(1999). 

\bibitem{cf:DemZei}
A.~Dembo and O.~Zeitouni, 
\textit{Large deviations techniques and applications},
Applications of Mathematics, \textbf{38}, Springer--Verlag (1998).

\bibitem{cf:DS}
J.D.~Deuschel and D.W.~Stroock,
\textit{Large Deviations},
Academic Press, Series in Pure
and Applied Mathematics, \textbf{137} (1989).

\bibitem{cf:FFS}
R.~Fernandez, J.~Fr\"ohlich and A.D.~Sokal
\textit{Random walks, critical phenomena,  and triviality
in quantum field theory},
Texts and monographs in physics, Springer--Verlag
 (1992).


 
\bibitem{cf:Georgii}
H.--O.~Georgii,
\textit{Gibbs Measures and Phase Transitions},  Studies in Mathematics,
\textbf{9}, W.~de Gruyter ed.
 (1988).   




 \bibitem{cf:HP}
I.~Herbst and L.D.~Pitt,
\textit{Diffusion equation techniques in stochastic monotonicity and positive
correlations}, Probab. Theory Rel. Fields,
\textbf{87} (1991), 275--312.


\bibitem{cf:Lawler}
G.F.~Lawler,
\textit{Intersections of Random Walks},
Probability and its Applications, Birkh\"auser,
 (1991).
 
 \bibitem{cf:LM}
 J.L.~Lebowitz and C.~Maes,
\textit{ The effect of an external field on an interface, entropy repulsion},
 J. Stat. Phys.
\textbf{46} (1987),
 39--49.      


 \bibitem{cf:Preston}
C.J.~Preston,
\textit{A generalization of the FKG inequalities}, 
 Comm. Math. Phys. \textbf{36} (1974),
 233-241.   


         







 
         

\end{thebibliography}
\end{document}